\documentclass[a4paper]{scrartcl}
\usepackage{amsmath,amssymb,url,graphicx, verbatim}
\usepackage{varioref}
\usepackage{relsize}
\newtheorem{theorem}{Theorem}[section]
\newtheorem{lemma}[theorem]{Lemma}

\newtheorem{corollary}[theorem]{Corollary}

\DeclareMathOperator{\abs}{abs}
\numberwithin{equation}{section}
\newcommand{\extcite}[2]{\cite[#2]{#1}}

\begin{document}

\begin{center}
 \LARGE \textbf{STAR DISCREPANCY BOUNDS OF DOUBLE INFINITE MATRICES INDUCED BY LACUNARY SYSTEMS}
\end{center}

\vspace{3ex}

\begin{center}
 THOMAS L\"OBBE
\footnote{The results are part of the author's PhD thesis supported by IRTG 1132, University of Bielefeld}
\end{center}

\vspace{3ex}

\begin{abstract}
ABSTRACT. In 2001 Heinrich, Novak, Wasilkowski and Wo\'zniakowski proved that the inverse of the star discrepancy satisfies $n(d,\varepsilon)\leq c_{\abs}d \varepsilon^{-2}$ by showing that there exists a set of points in $[0,1)^d$ whose star-discrepancy is bounded by $c_{\abs}\sqrt{d/N}$.  This result was generalized by Aistleitner who showed that there exists a double infinite random matrix with elements in $[0,1)$ which partly are coordinates of elements of a Halton sequence and partly independent uniformly distributed random variables such that any $N\times d$-dimensional projection defines a set $\{x_1,\ldots,x_N\}\subset [0,1)^d$ with
\begin{equation*}
 D^*_N(x_1,\ldots,x_N)\leq c_{\abs}\sqrt{d/N}.
\end{equation*}
In this paper we consider a similar double infinite matrix where the elements instead of independent random  variables are taken from a certain multivariate lacunary sequence and prove that with high probability each projection defines a set of points which has up to some constant the same upper bound on its star-discrepancy but only needs a significantly lower number of digits to simulate.
\end{abstract}

\vspace{5ex}

\section{Introduction}
\subsection*{Discrepancy and Uniform Distribution}

A sequence of vectors $(x_n)_{n\geq 1}=(x_{n,1},\ldots,x_{n,d})_{n\geq 1}$ of real numbers in $[0,1)^d$ is called \textit{uniformly distributed modulo one} if
\begin{equation}
 \label{81}
 \lim_{N\to\infty}\frac{1}{N}\sum_{n=1}^N\mathbf{1}_{\mathcal{A}}(x_n)=\lambda(\mathcal{A})
\end{equation}
for any axis-parallel box $\mathcal{A}\subset [0,1)^d$ where $\mathbf{1}_{\mathcal{A}}$ denotes the indicator function on the set $\mathcal{A}$ and $\lambda$ denotes the Lebesgue-measure on $[0,1)^d$.
The discrepancy resp. the star discrepancy of the first $N$ elements of $(x_n)_{n\geq 1}$ is defined by
\begin{equation}
 \label{65}
 \begin{aligned}
 D_N(x_1,\ldots,x_N)\:\:\: & = \:\:\: \sup_{\mathcal{A}\in\mathcal{B}}\left|\frac{1}{N}\sum_{n=1}^N\mathbf{1}_{\mathcal{A}}(x_n)-\lambda(\mathcal{A})\right|,\\
 D^*_N(x_1,\ldots,x_N)\:\:\: & =\:\:\: \sup_{\mathcal{A}\in\mathcal{B}^*}\left|\frac{1}{N}\sum_{n=1}^N\mathbf{1}_{\mathcal{A}}(x_n)-\lambda(\mathcal{A})\right|,
 \end{aligned}
\end{equation}
where $\mathcal{B}$ denotes the set of all axis-parallel boxes $\mathcal{A}=\prod_{i=1}^d[\alpha_{i},\beta_{i})\subset [0,1)^d$ and furthermore
$\mathcal{B}^*$ denotes the set of all axis-parallel boxes $\mathcal{A}=\prod_{i=1}^d[0,\beta_{i})\subset [0,1)^d$ with one corner in $0$. 
It is well-known that (\ref{81}) is equivalent to $D_N(x_1,\ldots,x_N)\to 0$ resp. $D^*_N(x_1,\ldots,x_N)\to 0$ for $N\to\infty$. By a classical result of Weyl \cite{W16} it is known that for any increasing sequence $(M_n)_{n\geq 1}$ of positive integers the sequence $(\langle M_nx\rangle)_{n\geq 1}$, where $\langle\cdot\rangle$ denotes the fractional part, is uniformly distributed modulo one for almost all $x\in [0,1)$. This result naturally extends to the multidimensional case. Sequences with vanishing star-discrepancy have applications in the theory of numerical integration. The connection is established by the Koksma-Hlawka inequality (see \cite{DT97}) which states that for any sequence of vectors $(x_n)_{n\geq 1}\subset [0,1)^d$ we have
\begin{equation}
 \label{217}
 \left|\frac{1}{N}\sum_{n=1}^Nf(x_n)-\int_{[0,1)^d}f(x)\,dx\right|\leq D_N^*(x_1,\ldots,x_N)\cdot V_{HK}(f)
\end{equation}
for any function $f$ on $[0,1)^d$ where $V_{HK}$ denotes the total variation in the sense of Hardy and Krause. Thus the integral can be approximated by the mean of the values which some points have under $f$ where the approximation error is given by the total variation of $f$ and the star-discrepancy of the points.
Therefore we are not only interested in sequences such that the star-discrepancy tends to 0, but also in the speed of convergence.


\subsection*{Low-discrepancy sequences}

Now we introduce Halton sequences which extend the definition of Van der Corput sequences to the multidimensional case. For an integer $d\geq 1$ let $(p_i)_{1\leq i\leq d}$ be a system of $d$ pairwise coprime integers. Then for any integer $n$ and $i\in\{1,\ldots,d\}$ let the $p_i$-adic decomposition of $n$ be given by
\begin{equation*}
 n=\sum_{j=0}^{\infty}\alpha(j,i)p_i^j
\end{equation*}
with $\alpha(j,i)\in\{0,\ldots,p_i-1\}$ for all $j\geq 0$ where $|\{j\in\mathbb{N}\cup\{0\}:\alpha(j,i)\neq 0\}|<\infty$ for all $n\in\mathbb{N}$ and $i\in\{1,\ldots,d\}$. For $n\in\mathbb{N}$ and $i\in\{1,\ldots,d\}$ define
\begin{equation*}
 x_{n,i}=\sum_{j=0}^{\infty}\alpha(j,i)p_i^{-j-1}.
\end{equation*}
The sequence $(x_n)_{n\geq 1}=(x_{n,1},\ldots,x_{n,d})_{n\geq 1}$ is called a Halton sequence in base $(p_1,\ldots,p_d)$. The discrepancy of a Halton sequence satisfies
\begin{equation}
 \label{218}
 D_N(x_1,\ldots,x_N)\leq C_d\frac{\log(N)^d}{N}
\end{equation}
with some constant $C_d>0$ which depends on $d$ (see \cite{H60}). Observe that the one-dimensional projection can be represented as the orbit of a von Neumann-Kakutani transformation. By using a randomly chosen starting point for this transformation Wang and Hickernell \cite{WH00} introduced the so-called randomized Halton sequences. Sequences with a discrepancy satisfying (\ref{218}) are called low-discrepancy sequences. Numerical integration using deterministic low-discrepancy sequences is called Quasi-Monte Carlo (QMC) integration in contrast to classical Monte Carlo integration which uses independent randomly chosen points. Many examples of deterministic low-discrepancy sequences can be found in the books of Dick and Pillichshammer \cite{DP10} and also Niederreiter \cite{N92}. A lower bound on the discrepancy was given by Roth \cite{R80} who proved
\begin{equation*}
 D_N(x_1,\ldots,x_N)\geq C_d\frac{\log(N)^{d/2}}{N}
\end{equation*}
for infinitely many $N$, some constant $C_d>0$ depending only on $d$ and any sequence of points $(x_n)_{n\geq 1}$.\\
Although low-discrepancy sequences have best known asymptotic bounds there are difficulties in applying them in practice. There are many applications which demand evaluation of high-dimensional integrals.
The upper bound on the right-hand side of (\ref{218}) only is vanishing if $N\geq e^d$ and thus such an upper bound is not feasible for high-dimensional integration in practice. There are some particular low-discrepancy sequences which provide good results in some special applications. For example, Atanassov \cite{A04} modified the definition 
of a Halton sequence obtaining a constant $C_d$ on the right-hand side of (\ref{218}) vanishing exponentially in $d$. But in general the situation is dissatisfying. Therefore randomized Quasi-Monte Carlo methods were introduced which try to combine the advantages of Quasi-Monte Carlo methods and classical Monte Carlo methods. Observe that the latter ones provide error bounds which are independent of the dimension while the former ones provide good asymptotic error 
bounds. Randomized Halton sequences are one example. For further example see the book of Lemieux \cite{L09} and the references therein.

\subsection*{Inverse of the star-discrepancy}

Since low-discrepancy sequence only give good error bounds if the number of points is large in comparison with the dimension, one could ask about sequences which have small discrepancy in the special case of a ``small'' number of sample points in comparison with the dimension. This led to the introduction of the ``inverse of the star-discrepancy''
\begin{equation*}
 n(d,\varepsilon)=\min\{N\in\mathbb{N}:\exists x_1,\ldots,x_N\in [0,1)^d,D^*_N(x_1,\ldots,x_N)\leq\varepsilon\}
\end{equation*}
which states the smallest number of points in $[0,1)^d$ having the upper bound $\varepsilon$ on the star-discrepancy. Heinrich, Novak, Wasilkowski and Wo\'zniakowski \cite{HNWW01} showed
\begin{equation}
 n(d,\varepsilon)=\mathcal{O}(d\varepsilon^{-2})
\end{equation}
with some implied constant which is independent of $d$ and $\varepsilon$. Thus there is a sequence of points in $[0,1)^d$ with
\begin{equation}
 \label{219}
 D^*_N(x_1,\ldots,x_N)\leq C\frac{\sqrt{d}}{\sqrt{N}}
\end{equation}
which for small $N$ compared with $d$ gives a better bound than (\ref{218}). Furthermore Hinrichs \cite{H04} proved
\begin{equation}
 n(d,\varepsilon)=\Omega(d\varepsilon^{-1}).
\end{equation}
Thus the dependence of $d$ in (\ref{219}) is optimal, only the precise order of $\varepsilon$ is unknown. In applications it often is desirable to have a sequence which is extendable not only in the number of points but also in dimension. Therefore 
Dick \cite{D07} proved that there exists a double infinite matrix $(x_{n,i})_{n\geq 1,i\geq 1}$ with numbers $x_{n,i}\in [0,1)$ such that for any pair of natural numbers $N,d\geq 1$ the projection $(x_{n,i})_{1\leq n\leq N, 1\leq i\leq d}$ defines an $N$-element sequence of points 
\begin{equation*}
\left\{(x_{1,1},\ldots,x_{1,d}),\ldots,(x_{N,1},\ldots,x_{N,d})\right\}\subset [0,1)^d
\end{equation*}
with star-discrepancy
\begin{equation}
 D^*_N(x_1,\ldots,x_N)\leq C\frac{\sqrt{d\log(N)}}{\sqrt{N}}
\end{equation}
for some absolute constant $C>0$ independent of $d$. Observe that the logarithmic term is due to the fact that Dick actually proved that any matrix generated by independent uniformly distributed random variables satisfies the upper bound with positive probability. The result was later improved by Doerr, Gnewuch, Kritzer and Pillichshammer \cite{DGKP08} who showed
\begin{equation}
 D^*_N(x_1,\ldots,x_N)\leq C\frac{\sqrt{d\log(1+N/d)}}{\sqrt{N}}.
\end{equation}
Aistleitner and Weimar \cite{AW} later obtained
\begin{equation}
 D^*_N(x_1,\ldots,x_N)\leq \frac{\sqrt{C_1d+C_2\log(\log(N))}}{\sqrt{N}}
\end{equation}
which is the best possible result because of the Chung-Smirnov Law of the Iterated Logarithm.

To avoid the iterated logarithm term in the upper bound hybrid sequences which are partly constructed by random numbers and partly by elements of a low-discrepancy sequence were introduced. Aistleitner \cite{A13} constructed a matrix where for large $n$ compared to $i$ the entries $x_{n,i}$ are taken from a Halton sequence while for small $n$ compared to $i$ they are randomly chosen. He proved that there exists a matrix which satisfies (\ref{219}) uniformly in $N$ and $d$.

\subsection*{Lacunary sequences}

Let $(A_n)_{n\geq 1}$ be a sequences of integer-valued $d\times d$-matrices and set $M_n=A_n\cdots A_1$ for all $n\geq 1$. Furthermore let $f:\mathbb{R}^d\to\mathbb{R}$ be a bounded periodic function of mean zero which is of bounded total variation in the sense of Hardy and Krause. By a result of Conze, Le Borgne, Roger \cite{CBR12} the system $(f(M_nx))_{n\geq 1}$ satisfies the Central Limit Theorem if
\begin{equation}
\label{17}
||M_{n+k}^Tj||_{\infty}\geq q^k||M_n^T||_{\infty}
\end{equation}
for all $j\in\mathbb{Z}^d\backslash\{0\}$, $n\in\mathbb{N}$, $k\geq \log_q(||j||_{\infty})$ and some absolute constant $q>1$
In general sequences which satisfy (\ref{17}) are called multivariate lacunary sequences.
For $x\in\mathbb{R}^d$ denote the vector which entries are the fractional parts of the entries of $x$ by $\langle x\rangle$.
Then the centered indicator functions $f_{\mathcal{A}}(\cdot)=\mathbb{1}_{\mathcal{A}}(\cdot)-\lambda{\mathcal{A}}$ for a suitable system of boxes $\mathcal{A}$ by which the discrepancy resp. the star-discrepancy are defined are typical example for periodic functions for which the Central Limit Theorem holds. 
Thus the lacunary system $(f_{\mathcal{A}}(M_nx))_{n\geq 1}$ shows a behaviour typical for independent, identically distributed random variables.
In applications it is more reasonable to use a point set defined by a lacunary system instead of set of independent random points. Both sequences have a similar probabilistic behaviour but the computational cost for simulating a suitable lacunary system is significantly smaller.

\subsection*{Main result}

We consider a similar constructed double infinite matrix. We define a double infinite matrix $(x_{n,i})_{n,i\geq 1}$ where $(x_{1,i})_{i\geq 1}$ forms a family of independent uniformly in $[0,1)$ distributed random variables. While for large $n$ compared to $i$ we define $x_{n+1,i}$ by taking elements of randomized Halton sequences, for small $n$ compared to $i$ we take fractional parts of a lacunary sequences instead of independent random numbers $x_{n,i}$. 

The practical purpose of having points defined by such a lacunary sequence instead of independent random points is reducing the number of digits which are necessary to simulate those points. To simulate $N$ random points in $[0,1)^d$ with a precision of $H$ digits requires a simulation of $dHN$ digits while by using points from such a lacunary sequence this number may be reduced to $\mathcal{O}(dH+d\log(d)N)$.





Before stating the main theorem we review the definition of a randomized Halton sequence.
Let $x=(x_1,\ldots,x_d)\in [0,1)^d$. For $i\in\{1,\ldots,d\}$ and some integer $p_i\geq 2$ we define the $p_i$-adic decomposition of $x_i$ by $x_i=\sum_{j=0}^{\infty}\alpha(j,i)p_i^{-j-1}$ where $\alpha(j,i)\in\{0,\ldots,p_i-1\}$ for all $i\in\{1,\ldots,d\}$ and $j\geq 0$. Now set
\begin{equation*}
 T_{p_i}(x_i)=\frac{\alpha(m,i)+1}{p_i^{m+1}}+\sum_{j>m}\frac{\alpha(j,i)}{p_i^{j+1}}
\end{equation*}
where $m=\min\{j:\alpha(j,i)\neq p_i-1\}$. Furthermore for a collection of pairwise coprime odd integers $p=(p_1,\ldots,p_d)$ set $T_p(x)=(T_{p_1}(x_1),\ldots,T_{p_d}(x_d))$. Observe that for $x_0=0$ the sequence $(x_n)_{n\geq 1}$ with $x_n=T_p(x_{n-1})$ for $n\geq 1$ defines a Halton sequence. Therefore for some uniformly distributed $x_0\in [0,1)^d$ and pairwise coprime odd integers $p_1,\ldots,p_d$ we call this sequence a randomized Halton sequence.

\begin{theorem}
 \label{164}
 Let $(x_{1,i})_{i\geq 1}$ be a sequence of independent random variables which are uniformly distributed in $[0,1)$ and let $(p_i)_{i\geq 1}$ be the sequence of all odd prime numbers. For all integers $n\geq 1$ and $i\geq 1$ define
 \begin{equation}
  \label{165}
  x_{n+1,i}=
  \begin{cases}
   T_{p_i}(x_{n,i}), & \textnormal{if }i=1\textnormal{ or }2^{12\cdot 2^i}<n-1,\\
   \langle 2^{\lceil\log_2(i)\rceil+1}x_{n,i}\rangle, & \textnormal{if }i\geq 2\textnormal{ and }2^{12\cdot 2^i}\geq n-1.
  \end{cases}
 \end{equation}
Then for any $\varepsilon>0$ the probability, that for any integers $N\geq 1$ and $d\geq 1$ the set of points $P=\{(x_{1,1},\ldots,x_{1,d}),\ldots,(x_{N,1},\ldots,x_{N,d})\}\subset [0,1)^d$ satisfies
\begin{equation}
 \label{166}
 D^*_N(P)\leq (2576+357\log(\varepsilon^{-1}))\frac{\sqrt{d}}{\sqrt{N}},
\end{equation}
is at least $1-\varepsilon$.
\end{theorem}

\vspace{5ex}

\section{Preliminaries}
\label{S41}

\begin{lemma}[Maximal Bernstein inequality, \extcite{EM96}{Lemma 2.2}]
 \label{191}
 For some integer $N\geq 1$ let $Z_1,\ldots,Z_N$ be a sequence of i.i.d. random variables with mean zero and variance $\sigma^2>0$ such that $|Z_1|\leq 1$. Then for any $t>0$ we have
 \begin{equation}
  \label{195}
   \mathbb{P}\left(\max_{M\in\{1,\ldots,N\}}\left|\sum_{n=1}^MZ_n\right|>t\right)\leq 2\exp\left(-\frac{t^2}{2N\sigma^2+2t/3}\right).
 \end{equation}
\end{lemma}
 
For integers $N\geq 1$ and $d\geq 1$ and an $N$-element set of $d$-dimensional points
\begin{equation*}
 \left\{(x_{1,1},\ldots,x_{1,d}),\ldots,(x_{N,1},\ldots,x_{N,d})\right\}
\end{equation*}
denote the star-discrepancy by $D^d_N(x_{n,i})$. Furthermore for an integer $0\leq M<N$ write $D_{M,N}^d(x_{n,i})$ for the star-discrepancy of the $N-M$-element point set
\begin{equation*}
 \left\{(x_{M+1,1},\ldots,x_{M+1,d}),\ldots,(x_{N,1},\ldots,x_{N,d})\right\}.
\end{equation*}

\begin{lemma}[\extcite{DP10}{Proposition 3.16}]
\label{168}
Let $0\leq M<N$ be integers. Then for points $y_1,\ldots,y_N\in [0,1)^d$ we have
\begin{equation*}
 D^d_N(y_1,\ldots,y_N)\leq\frac{MD^d_M(y_1,\ldots,y_M)}{N}+ \frac{(N-M)D^d_{M,N}(y_{M+1},\ldots,y_N)}{N}
\end{equation*}
and
\begin{equation*}
 D^d_{M,N}(y_{M+1},\ldots,y_N)\leq\frac{ND^d_N(y_1,\ldots,y_N)}{N-M}+ \frac{MD^d_M(y_{1},\ldots,y_M)}{N-M}.
\end{equation*}
\end{lemma}

Let $v,w\in [0,1)^d$. We write $v\leq w$ if $v_i\leq w_i$ for all $i\in\{1,\ldots,d\}$. For some $\delta>0$ a set $\Delta$ of elements in $[0,1)^d\times [0,1)^d$ is called a $\delta$-bracketing cover if for every $x\in [0,1)^d$ there exists $(v,w)\in\Delta$ with $v\leq x\leq w$ and $\lambda(\overline{[v,w)})\leq\delta$ for $\overline{[v,w)}=[0,w)\backslash [0,v)$.
The following Lemma gives an upper bound on the cardinality of a $\delta$-bracketing cover.

\begin{lemma}[\extcite{G08}{Theorem 1.15}]
 \label{176}
 For any $d\geq 1$ and $\delta>0$ there exists some $\delta$-bracketing cover $\Delta$ with
 \begin{equation*}
  |\Delta|\leq \frac{1}{2}(2e)^d(\delta^{-1}+1)^d.
 \end{equation*}
\end{lemma}

\begin{corollary}
 \label{177}
 For any integers $d\geq 1$ and $h\geq1$ there exists a $2^{-h}$-bracketing cover $\Delta$ with
 \begin{equation*}
  |\Delta|\leq \frac{1}{2}(2e)^d(2^{h+2}+1)^d
 \end{equation*}
 such that for any $(v,w)\in\Delta$ and any $i\in \{1,\ldots,d\}$ we have
 \begin{eqnarray*}
  v_i & = & 2^{-(\lceil\log_2(i)\rceil +1)(h+1)}a_i,\\
  w_i & = & 2^{-(\lceil\log_2(i)\rceil +1)(h+2)}b_i
 \end{eqnarray*}
 for some integers $a_i\in\{0,1,\ldots,2^{(\lceil\log_2(i)\rceil +1)(h+1)}\}$ and $b_i\in\{0,1,\ldots,2^{(\lceil\log_2(i)\rceil +1)(h+2)}\}$.
\end{corollary}

\textit{Proof.} Let $\Delta$ be some $2^{-(h+2)}$-bracketing cover of $[0,1)^d$. By Lemma \ref{176} we have
\begin{equation*}
 |\Delta|\leq \frac{1}{2}(2e)^d(2^{(h+2)}+1)^d.
\end{equation*}
For $(v,w)\in\Delta$ and $i\in\{1,\ldots,d\}$ define
\begin{eqnarray*}
 y_{v,i} & = & \max\left\{2^{-(\lceil\log_2(i)\rceil +1)(h+1)}a_i\leq v_i:a_i\in\mathbb{Z}\right\},\\
 z_{w,i} & = & \min\left\{2^{-(\lceil\log_2(i)\rceil +1)(h+2)}b_i\geq w_i:b_i\in\mathbb{Z}\right\}.
\end{eqnarray*}
For $y_v=(y_{v,i})_{i\in\{1,\ldots,d\}}\in [0,1)^d$ we obtain
\begin{equation*}
 \lambda(\overline{[y_v,v)})\leq \sum_{i=1}^d2^{-(\lceil\log_2(i)\rceil +1)(h+1)}\leq 2^{-(h+1)}\sum_{i=1}^di^{-(h+1)}\leq 2^{-(h+1)}.
\end{equation*}
Analogously for $z_w=(z_{w,i})_{i\in\{1,\ldots,d\}}\in [0,1)^d$ we have
\begin{equation*}
 \lambda(\overline{[z,z_w)})\leq 2^{-(h+2)}.
\end{equation*}
Thus we get
\begin{equation*}
 \lambda(\overline{[y_v,z_w)})\leq \lambda(\overline{[y_v,v)})+\lambda(\overline{[v,w)})+\lambda(\overline{[w,z_w)})\leq 2^{-h}.
\end{equation*}
Set $\tilde{\Delta}=\{(y_v,z_w):(v,w)\in\Delta\}$. Since $\Delta$ is a $2^{-(h+2)}$-bracketing cover for any $x\in[0,1)^d$ there exists $(v,w)\in\Delta$ and $(y_v,z_w)\in\tilde{\Delta}$ with $y_v\leq v\leq x\leq w\leq z_w$. Therefore $\tilde{\Delta}$ is a $2^{-h}$-bracketing cover and the conclusion of the proof follows by $|\tilde{\Delta}|\leq |\Delta|$.

\vspace{5ex}

\section{Randomized Halton sequences}
\label{S42}

Note that we assume that the integers $p_1,\ldots, p_d$ are odd since we later need sequences such that not only $(x_n)_{n\geq 1}$ is a low-discrepancy sequence but also subsequences $(x_{n_l})_{l\geq 1}$ where the elements $n_l$ belong to one particular modulo class with modulo $2^{\kappa}$ for some integer $\kappa$ have sufficiently small discrepancy. This shall be ensured by the following

\begin{lemma}
 \label{170}
 For some integer $d\geq 2$ let $(x_n)_{n\geq 1}$ be a randomized Halton sequence in $[0,1)^d$ constructed by the first $d$ odd primes. Let $N\geq 2^{12\cdot 2^{d}}$ be the number of points. For some integers $1\leq \kappa\leq \log_2(8\log_2(N))$ and $\gamma\in\{0,\ldots,2^{\kappa}-1\}$ set 
 \begin{equation*}
  N_{\kappa,\gamma}=\{n:n\in\{1,\ldots,N\},n\equiv\gamma \pmod{2^{\kappa}}\}
 \end{equation*} 
 and define $P_{N,\kappa,\gamma}=\{x_n:n\in N_{\kappa,\gamma}\}$. Then the star-discrepancy of $P_{N,\kappa,\gamma}$ satisfies
 \begin{equation}
 \label{149} 
  D_{N,\kappa,\gamma}^*(\{x_n:n\in N_{\kappa,\gamma}\})=D_{|N_{\kappa,\gamma}|}^*(P_{N,\kappa,\gamma})\leq \frac{\sqrt{d}}{\sqrt{|N_{\kappa,\gamma}|}}.
 \end{equation}

\end{lemma}

\textit{Proof.} The proof of this Lemma which is an application of the Chinese Remainder Theorem is mainly based on the proofs of \extcite{N92}{Theorem 3.6} and \extcite{A13}{Corollary 1}. For some $x_0$ let the $p_i$-adic decomposition of $x_{n,i}=(T_{p_i}^{n}(x_0))_i$ be given by 
\begin{equation*}
 x_{n,i}=\sum_{j=0}^{\infty}\alpha(j,i,x_0,n)p_i^{-j-1}
\end{equation*} 
for suitable integers $\alpha(j,i,x_0,n)\in\{0,1,\ldots,p_i-1\}$. Observe that for any $i\in\{1,\ldots,d\}$ there exists at most one $N_i\in\{1,\ldots,N\}$ such that
\begin{multline*}
 \alpha(j,i,x_0,1)=\alpha(j,i,x_0,2)=\cdots=\alpha(j,i,x_0,N_i)\\
 \neq \alpha(j,i,x_0,N_i+1)=\cdots =\alpha(j,i,x_0,N)
\end{multline*}
for all $j\geq\lceil \log_{p_i}(N)\rceil$. Let $\pi:\{1,\ldots,d\}\to\{1,\ldots,d\}$ be a permutation satisfying $N_{\pi(i)}\geq N_{\pi(j)}$ if $i\geq j$ and set $\pi(0)=0$ and $\pi(d+1)=N$. Therefore there exist constants $g_{m,i}$ for any $m\in\{1,\ldots,d+1\}$ and $i\in\{1,\ldots,d\}$ such that
\begin{equation}
 \label{156}
 \sum_{j=\lceil \log_{p_i}(N)\rceil}^{\infty}\alpha(j,i,x_0,n){p_i}^{-j-1} =g_{m,i}
\end{equation}
for all $n\in\{N_{\pi(m-1)}+1,\ldots,N_{\pi(m)}\}$. Now fix some set $\mathcal{N}=\{N_{\pi(m-1)}+1,\ldots,N_{\pi(m)}\}$ and define
\begin{equation*}
 \Psi_i(x_n)=\Psi_i\left(\sum_{j=0}^{\lceil \log_{p_i}(N)\rceil-1}\alpha(j,i,x_0,n){p_i}^{-j-1}+g_{m,i}\right)=\sum_{j=0}^{\lceil \log_{p_i}(N)\rceil-1}\alpha(j,i,x_0,n){p_i}^j.
\end{equation*}
Set $n_l=2^{\kappa}(l-1)+\gamma$ for any integer $l\geq 1$ and some $\gamma\in\{0,\ldots,2^{\kappa}-1\}$.
By definition of $T_p$ it is easy to see that for $n,n+1\in\mathcal{N}$ and $i\in\{1,\ldots,d\}$ we have $\Psi_i(x_{n+1})=\Psi_i(x_{n})+1$. Thus we obtain $\Psi_i(x_{n_{l+1}})=\Psi_i(x_{n_l})+2^{\kappa}$ for $n_{l+1},n_l\in\mathcal{N}$. We now shall show the following version of the Chinese Remainder Theorem:\\
Let $\beta_1,\ldots,\beta_d$ and $s_1,\ldots,s_d$ be positive numbers, then there exists an integer $\beta$ such that any solution of
\begin{equation}
 \label{151}
 \begin{aligned}
  \Psi_1(x_{n_l}) \:\:\: & \equiv \:\:\: \beta_1 \quad \pmod {p_1^{s_1}}\\
  & \:\: \vdots \\
  \Psi_d(x_{n_l}) \:\:\: & \equiv \:\:\: \beta_d \quad \pmod {p_d^{s_d}}
 \end{aligned}
\end{equation}
satisfies 
\begin{equation}
 \label{152}
 \quad\quad\quad\quad\quad\quad\quad\:\:\, l \:\:\: \equiv \:\:\: \beta \quad \:\:\pmod {p_1^{s_1}p_2^{s_2}\cdots p_d^{s_d}}.
\end{equation}
Observe that we have $\Psi_i(x_{n_{l'}})\equiv \Psi_i(x_{n_l}) \,\, \pmod {p_i^{s_i}}$ only if $2^{\kappa}(l'-l)\equiv 0 \,\, \pmod {p_i^{s_i}}$. Since $2^{\kappa}$ and $p_i^{s_i}$ are coprime we have $(l'-l)|p_i^{s_i}$. Now for any $i\in\{1,\ldots,d\}$ define the map $\Xi_i:\mathbb{Z}\to\mathbb{Z}/p_i^{s_i}\mathbb{Z}$ with $\Xi_i(l)=\Psi_i(x_{n_l})+p_i^{s_i}\mathbb{Z}$. Observe that $\Xi_i$ is periodic, i.e. $\Xi_i(l)=\Xi_i(l+z p_i^{s_i})$ for any integer $z$, and $\Xi_i|_{\{1,\ldots,p_i^{s_i}\}}$ is bijective, i.e. $\Xi_i(l)\neq \Xi_i(l')$ for $l,l'\in\{1,\ldots,p_i^{s_i}\}$ with $l\neq l'$. Since the system (\ref{151}) only has a solution if $\Xi_i(l)= \beta_i +p_i^{s_i}$ for all $i\in\{1,\ldots,d\}$ we see that there are integers $\alpha_1,\ldots,\alpha_d$ such any solution satisfies $l\equiv \alpha_i \,\, \pmod {p_i^{s_i}}$ for $i\in\{1,\ldots,d\}$.
By classical Chinese Remainder Theorem there exists some integer $\beta$ such for any solution $l$ we conclude (\ref{152}).
Thus among $\prod_{i=1}^dp_i^{s_i}$ consecutive numbers of the sequence $(n_l)_{n\geq 1}$ there is exactly one $l$ such that $\Psi_i(x_{n_l})=a_i$ for any collection of numbers $a_i\in\{0,\ldots,p_i^{s_i}-1\}$ with $i\in\{1,\ldots,d\}$. 
Let $B$ be any box of the form
\begin{equation*}
B=\prod_{i=1}^d[a_ip_i^{-s_i},(a_i+1)p_i^{-s_i})
\end{equation*}
with $a_i\in\{0,\ldots,p_i^{s_i}-1\}$ and $i\in\{1,\ldots,d\}$.
Observe that $x\in B$ if for $i\in\{1,\ldots,d\}$ the first $s_i$ digits in the $p_i$-adic decomposition of $x_i$ are uniquely defined, i.e. we have $\sum_{j=0}^{s_i-1}\alpha(j,i)p_i^{-j-1}=a_i$ for all 
\begin{equation*}
 x_i=\sum_{j=0}^{\infty}\alpha(j,i)p_i^{-j-1}\in B_i=[a_ip_i^{-s_i},(a_i+1)p_i^{-s_i}).
\end{equation*}
 By definition of the $\Psi_i$ this is equivalent to $\Psi_i(x)\equiv \sum_{j=0}^{s_i-1}\alpha(j,i)p_i^{j} \,\, \pmod {p_i^{s_i}}$ for all $i\in\{1,\ldots,d\}$. Thus there is exactly one $l$ with $x_{n_l}\in B$ among $\prod_{i=1}^dp_i^{s_i}$ consecutive numbers of the sequence $(x_{n_l})_{l\geq 1}$. Therefore we obtain
\begin{equation}
 \label{153}
 \left|\left\{l:x_{n_l}\in B,l\in\left\{L+1,\ldots,L+\prod_{i=1}^dp_i^{s_i}\right\}\subset\mathcal{N}\right\}\right|=1.
\end{equation}

Now for $i\in\{1,\ldots,d\}$ and integers $r_i\geq 0$ let
\begin{equation*}
 \mathcal{C}_i(r_i)=\left\{[0,c_ip_i^{-r_i}):c_i\in\{0,\ldots, p_i^{r_i}-1\}\right\}
\end{equation*}
be a family of intervals. Furthermore set
\begin{equation*}
 \mathcal{A}_i(r_i)=\left\{[a_ip_i^{-s_i},(a_i+1)p_i^{-s_i}):a_i\in\{0,1\ldots,p_i^{s_i}-1\},s_i\in\{0,\ldots,r_i\}\right\}.
\end{equation*}
For integers $r_1,\ldots,r_d\geq 0$ let
\begin{equation*}
 \mathcal{B}(r_1,\ldots,r_d)=\left\{B=\prod_{i=1}^dB_i:B_i\in\mathcal{C}_i(r_i)\cup\mathcal{A}_{i}(r_i)\textnormal{ for any }i\in\{1,\ldots,d\}\right\}
\end{equation*}
be a collection of boxes.

For any box $B\subseteq [0,1)^d$ and any set $\mathcal{M}$ of positive integers set
\begin{equation*}
 \mathcal{D}_{\mathcal{M}}(B)=\left|\sum_{n\in \mathcal{M}}(\mathbf{1}_{B}(x_n)-\lambda(B))\right|.
\end{equation*}
Furthermore for integers $\kappa,\gamma\geq 0$ with $\gamma\in\{0,\ldots,2^{\kappa}-1\}$ and any set $\mathcal{N}$ as defined above let
\begin{equation*}
 \mathcal{N}_{\kappa,\gamma}=\left\{n\in\mathcal{N}:n\equiv\gamma \,\, \pmod {2^{\kappa}}\right\}.
\end{equation*}
Now we shall show that for any set of integers $r_1,\ldots,r_d\geq 0$ and any box 
\begin{equation*}
 B=\prod_{i=1}^dB_i\in \mathcal{B}(r_1,\ldots,r_d)
\end{equation*}
we have
\begin{equation}
 \label{154}
 \mathcal{D}_{\mathcal{N}_{\kappa,\gamma}}(B)\leq \prod_{\substack{i\in\{1,\ldots,d\},\\ B_i\notin \mathcal{A}_i(r_i)}}\left(\frac{p_i-1}{2}r_i+1\right).
\end{equation}
We are going to prove this inequality by induction on the number $k$ of indices $i$ such that $B_i\notin\mathcal{A}_i(r_i)$. Thus we first assume $k=0$. We have $B=\prod_{i=1}^d[a_ip_i^{-s_i},(a_i+1)p_i^{-s_i})$ for suitable integers $s_1,\ldots,s_d$ and $a_1,\ldots,a_d$. By (\ref{153}) we obtain
\begin{equation*}
 \left\lfloor|\mathcal{N}_{\kappa,\gamma}|\prod_{i=1}^dp_i^{-s_i}\right\rfloor\leq \sum_{n\in \mathcal{N}_{\kappa,\gamma}}\mathbf{1}_B(x_n)\leq \left\lceil|\mathcal{N}_{\kappa,\gamma}|\prod_{i=1}^dp_i^{-s_i}\right\rceil.
\end{equation*}
Since $\sum_{n\in \mathcal{N}_{\kappa,\gamma}}\lambda(B)=|\mathcal{N}_{\kappa,\gamma}|\prod_{i=1}^dp_i^{-s_i}$ we conclude $\mathcal{D}_{\mathcal{N}_{\kappa,\gamma}}(B)\leq 1$ for $k=0$.
Now assume that (\ref{154}) has been proved for $|\{i:B_i\notin \mathcal{A}_i(r_i)\}|=k-1$. Consider some box $B\in\mathcal{B}(r_1,\ldots,r_d)$ with $|\{i:B_i\notin \mathcal{A}_i(r_i)\}|=k$. Without loss of generality we may assume $B_i\notin\mathcal{A}_i(r_i)$ for $i\in\{1,\ldots,k\}$ and $B_i\in\mathcal{A}_i(r_i)$ for $i\in\{k+1,\ldots,d\}$. Then we have
$B_k=[0,c_kp_k^{-r_k})$ for some integer $c_k$ with $c_k\in\{0,\ldots,p_k^{r_k}-1\}$. We get $c_kp_k^{-r_k}=\sum_{j=1}^{r_k}e_jp_k^{-j}$ for integers $e_j$ with $e_j\in\{0,\ldots,p_k-1\}$ for $1\leq j\leq r_k$. Therefore the interval $B_k$ can be decomposed into $e_1$ intervals of length $p_k^{-1}$, $e_2$ intervals of length $p_k^{-2}$ and so on. Set $e=\sum_{j=1}^{r_k}e_j$. Then
\begin{equation*}
 B_k=\bigcup_{t=1}^eE_t
\end{equation*}
for pairwise disjoint $E_t\in\mathcal{A}_k(r_k)$ with $t\in\{1,\ldots,e\}$. Thus we obtain
\begin{equation*}
 B=\bigcup_{t=1}^e(B_1\times\cdots\times B_{k-1}\times E_t\times B_{k+1}\times\cdots\times B_d).
\end{equation*}
By induction hypothesis we observe
\begin{equation}
 \label{155}
 \begin{aligned}
  \mathcal{D}_{\mathcal{N}_{\kappa,\gamma}}(B) \:\:\: & \leq \:\:\: \sum_{t=1}^e\mathcal{D}_{\mathcal{N}_{\kappa,\gamma}}\left(B_1\times\cdots\times B_{k-1}\times E_t\times B_{k+1}\times\cdots\times B_d\right)\\
  & \leq \:\:\: e\prod_{i=1}^{k-1}\left(\frac{p_i-1}{2}r_i+1\right).
 \end{aligned}
\end{equation}
Furthermore set $F=[c_kp_k^{-r_k},1)=[0,1)\backslash B_k$. We have
\begin{equation*}
 \begin{aligned}
  \mathcal{D}_{\mathcal{N}_{\kappa,\gamma}}(B) \:\:\: \leq & \:\:\: \mathcal{D}_{\mathcal{N}_{\kappa,\gamma}}\left(B_1\times\cdots\times B_{k-1}\times [0,1)\times B_{k+1}\times\cdots\times B_d\right)\\
  & \:\:\: +\mathcal{D}_{\mathcal{N}_{\kappa,\gamma}}\left(B_1\times\cdots\times B_{k-1}\times F\times B_{k+1}\times\cdots\times B_d\right).
 \end{aligned}
\end{equation*}
Thus we get
\begin{equation*}
 \begin{aligned}
  \mathcal{D}_{\mathcal{N}_{\kappa,\gamma}}(B) \:\:\: \leq & \:\:\: \prod_{i=1}^{k-1}\left(\frac{p_i-1}{2}r_i+1\right)\\
  \:\:\: & \:\:\: +\mathcal{D}_{\mathcal{N}_{\kappa,\gamma}}\left(B_1\times\cdots\times B_{k-1}\times F\times B_{k+1}\times\cdots\times B_d\right).
 \end{aligned}
\end{equation*}
Observe that $F$ can be decomposed into $1+\sum_{j=1}^{r_k}p_k-1-e_j=(p_k-1)r_k-e+1$ intervals in $\mathcal{A}_k(r_k)$. Thus we get
\begin{equation*}
 \mathcal{D}_{\mathcal{N}_{\kappa,\gamma}}(B)\leq ((p_k-1)r_k-e+2)\prod_{i=1}^{k-1}\left(\frac{p_i-1}{2}r_i+1\right).
\end{equation*}
By (\ref{155}) we have
\begin{equation*}
 \mathcal{D}_{\mathcal{N}_{\kappa,\gamma}}(B)\leq \frac{e+(p_k-1)r_k-e+2}{2}\prod_{i=1}^{k-1}\left(\frac{p_i-1}{2}r_i+1\right).
\end{equation*}
Hence (\ref{154}) is proved for any $k$.
Now let $J=\prod_{i=1}^d[0,v_i)\subset [0,1)^d$ be some arbitrary box. For any $i\in\{1,\ldots,d\}$ set $r_i=\lceil \log_{p_i}(N)\rceil$ and furthermore let $c_i$ be the integer such that $c_ip_i^{-r_i}\leq v_i<(c_i+1)p_i^{-r_i}$. Take some 
\begin{equation*}
 \mathcal{N}_{\kappa,\gamma,m}=\{n\in\{N_{\pi(m-1)}+1,\ldots,N_{\pi(m)}\},n\equiv \gamma \,\, \pmod {2^{\kappa}}\}.
\end{equation*}
By definition of $\mathcal{N}_{\kappa,\gamma,m}$ for any $i\in\{1,\ldots,d\}$ we have
$x_{n,i}=z_{n,i}p_i^{-r_i}+g_{m,i}$ for some integer $z_{n,i}$ depending on $n\in\mathcal{N}_{\kappa,\gamma,m}$ and $0\leq g_{m,i}<p_i^{-r_i}$ independent of $n$. For $v_i-c_ip_i^{-r_i}\leq g_{m,i}$ set $v'_i=c_ip_i^{-r_i}$, otherwise set $v'_i=(c_i+1)p_i^{-r_i}$ and let $B_{m}=\prod_{i=1}^d[0,v'_i)$. It is easy to see that
\begin{equation*}
 \sum_{n\in\mathcal{N}_{\kappa,\gamma,m}}\mathbf{1}_{J}(x_n)=\sum_{n\in\mathcal{N}_{\kappa,\gamma,m}}\mathbf{1}_{B_m}(x_n).
\end{equation*}
Thus by (\ref{154}) we get
\begin{eqnarray*}
 \mathcal{D}_{N_{\kappa,\gamma}}(J) & \leq & \sum_{m=1}^{d+1}\mathcal{D}_{\mathcal{N}_{\kappa,\gamma,m}}(J)\\
 & \leq & \sum_{m=1}^{d+1}\left|\sum_{n\in\mathcal{N}_{\kappa,\gamma,m}}(\mathbf{1}_{J}(x_n)-\lambda(J))\right|\\
 & \leq & \sum_{m=1}^{d+1}\left|\sum_{n\in\mathcal{N}_{\kappa,\gamma,m}}(\mathbf{1}_{B_m}(x_n)-\lambda(B_m))\right|+\sum_{m=1}^{d+1}|\mathcal{N}_{\kappa,\gamma,m}|\cdot|\lambda(J)-\lambda(B_m)|\\
 & \leq & \sum_{m=1}^{d+1}\mathcal{D}_{\mathcal{N}_{\kappa,\gamma,m}}(B_m)+|N_{\kappa,\gamma}|\sum_{i=1}^dp_i^{-r_i}\\
 & \leq & (d+1)\prod_{i=1}^d\left(\frac{p_i-1}{2}r_i+1\right)+|N_{\kappa,\gamma}|\sum_{i=1}^dp_i^{-r_i}.      
\end{eqnarray*}
Since $r_i=\lceil\log_{p_i}(N)\rceil$ we observe
\begin{equation}
 \label{157}
 \begin{aligned}
  D_{|N_{\kappa,\gamma}|}^*(P_{N,\kappa,\gamma}) \:\:\: & = \:\:\: \sup_{J}\frac{\mathcal{D}_{N_{\kappa,\gamma}}(J)}{|N_{\kappa,\gamma}|}\\
  & \leq \:\:\: \frac{d}{N}+\frac{1}{|N_{\kappa,\gamma}|}\prod_{i=1}^d\frac{i+1}{i}\left(\frac{p_i-1}{2\log(p_i)}\log(N)+\frac{p_i+1}{2}\right).
 \end{aligned}
\end{equation}
Next we shall show
\begin{equation}
 \label{158}
 \frac{i+1}{i}\left(\frac{p_i-1}{2\log(p_i)}\log(N)+\frac{p_i+1}{2}\right)\leq (i+1)\log(N).
\end{equation}
for all $i\in\{1,\ldots,d\}$. This is easy to see for $i\leq 4$. It is well-known that for $i\geq 5$ we have $i\leq p_i\leq 1+7/4\cdot i\log(i)$ (see, e.g. \extcite{BS96}{Theorem 8.8.4}).
Therefore we get
\begin{eqnarray*}
 \frac{i+1}{i}\left(\frac{p_i-1}{2\log(p_i)}\log(N)+\frac{p_i+1}{2}\right) & \leq & \frac{i+1}{i}\left(\frac{7/4\cdot i\log(i)}{2\log(i)}\log(N)+2i\log(i)\right)\\
 & \leq & \frac{i+1}{i}\left(\frac{7}{8} i\log(N)+2i\log(i)\right)\\
 & \leq & (i+1)\log(N).
\end{eqnarray*}
Thus (\ref{158}) is proved. Together with (\ref{157}) we have
\begin{equation}
 \label{159}
 D_{|N_{\kappa,\gamma}|}^*(P_{N,\kappa,\gamma})\leq \frac{d}{N}+\frac{(d+1)!(\log(N))^d}{|N_{\kappa,\gamma}|}.
\end{equation}
It remains to show
\begin{equation}
 \label{160}
 \frac{\sqrt{d}}{\sqrt{|N_{\kappa,\gamma}|}}+\frac{(d+1)!(\log(N))^{d}}{\sqrt{d|N_{\kappa,\gamma}|}}\leq 1.
\end{equation}
Then the statement of the Lemma follows by (\ref{159}) and
\begin{eqnarray*}
 D_{|N_{\kappa,\gamma}|}^*(P_{N,\kappa,\gamma}) & \leq & \frac{d}{|N_{\kappa,\gamma}|}+\frac{(d+1)!(\log(N))^d}{|N_{\kappa,\gamma}|}\\
 & \leq & \frac{\sqrt{d}}{\sqrt{|N_{\kappa,\gamma}|}}\left(\frac{\sqrt{d}}{\sqrt{|N_{\kappa,\gamma}|}}+\frac{(d+1)!(\log(N))^d}{\sqrt{d|N_{\kappa,\gamma}|}}\right)\\
 & \leq & \frac{\sqrt{d}}{\sqrt{|N_{\kappa,\gamma}|}}.
\end{eqnarray*}
In order to show (\ref{160}) we estimate the second term and observe
\begin{equation}
 \label{161}
 \frac{(d+1)!(\log(N))^{d}}{\sqrt{d|N_{\kappa,\gamma}|}}\leq \frac{(d+1)!(\log(N))^{d}}{\sqrt{2^{-\kappa-1}dN}} \leq \frac{4}{\sqrt{\log(2)}}\frac{(d+1)!(\log(N))^{d+1/2}}{\sqrt{dN}}
\end{equation}
where we used $\kappa\leq \log_2(8\log_2(N))$ for the second inequality.
Since for any fixed $d$ the derivative of $(\log(N))^{d+1/2}/\sqrt{N}$ is negative for $N\geq e^{2(d+1/2)}$ it is enough to restrict ourselves to the case $N=2^{12\cdot2^{d}}>e^{2(d+1/2)}$. Therefore we first shall show
\begin{equation}
 \label{162}
 \frac{4}{\sqrt{\log(2)}}\frac{(d+1)!}{\sqrt{d}}\leq \frac{1}{2}(\log(N))^{d-1/2}
\end{equation}
resp. equivalently
\begin{equation}
 \label{163}
 \frac{8(d+1)!}{12^{d-1/2}(\log(2))^d\sqrt{d}}\leq 2^{d^2-d/2}.
\end{equation}
This shall be done by induction. It can easily be verified that (\ref{163}) is true for $d=2$. Thus we may assume that (\ref{163}) holds for some integer $d\geq 2$. We get
\begin{eqnarray*}
 \frac{8((d+1)+1)!}{12^{(d+1)-1/2}(\log(2))^{(d+1)}\sqrt{d+1}} & \leq & (d+2)\frac{8(d+1)!}{12^{d-1/2}(\log(2))^d\sqrt{d}}\\
 & \leq & 2^{2d+1/2}\frac{8(d+1)!}{12^{d-1/2}(\log(2))^d\sqrt{d}}\\
 & \leq & 2^{2d+1/2}\cdot 2^{d^2-d/2}\\
 & \leq & 2^{(d+1)^2-(d+1)/2}
\end{eqnarray*}
and therefore we have (\ref{162}) for any integer $d\geq 2$. For $d\geq 2$ we get
\begin{equation*}
 (12\cdot 2^d)^{2d}\leq 2^{2d^2+2\log_2(12)d}\leq 2^{6\cdot 2^d}\leq \sqrt{N}.
\end{equation*}
Hence $(\log(N))^{2d}/\sqrt{N}\leq 1$ immediately follows. With $\sqrt{d}/\sqrt{|N_{\kappa,\gamma}|}\leq 1/2$ and (\ref{162}) we observe (\ref{160}) which finally concludes the proof.

\vspace{5ex}

\section{Proof of Theorem \ref{164}}
\label{S43}

The proof of this Theorem is mainly based on \cite{A13}. For some integers $N\geq 1$ and $d\geq 1$ we simply write
\begin{equation*}
 D^d_N(x_{n,i})=D^d_N((x_{1,1},\ldots,x_{1,d}),\ldots,(x_{N,1},\ldots,x_{N,d})).
\end{equation*}

For all integers $m\geq 1$ and $d\geq 1$ we define
\begin{equation*}
 \mathcal{F}_{m,d,\varepsilon}=\left\{\max_{M\in\{2^m+1,\ldots,2^{m+1}\}}MD^d_M(x_{n,i})\geq C_{m,d,\varepsilon}\sqrt{d}\sqrt{2^{m+1}}\right\}
\end{equation*}
with
\begin{equation*}
 C_{m,d,\varepsilon}=
 \begin{cases}
  1819+252\log(\varepsilon^{-1}), & \textnormal{if }12\cdot 2^d>m,\\
  1821+252\log(\varepsilon^{-1}), & \textnormal{if }12\cdot 2^d\leq m.
 \end{cases}
\end{equation*}
We shall show that
\begin{equation}
 \label{167}
 \mathbb{P}\left(\bigcup_{d\geq 1}\bigcup_{m\geq 1}\mathcal{F}_{m,d,\varepsilon}\right)\leq \varepsilon.
\end{equation}
Therefore on the complement of $\cup_{d\geq 1}\cup_{m\geq 1}\mathcal{F}_{m,d,\varepsilon}$ which has measure bounded from below by $1-\varepsilon$ for any integer $N\geq 1$ and $d\geq 1$ we have
\begin{equation*}
 ND^d_N(x_{n,i})\leq (1821+252\log(\varepsilon^{-1}))\sqrt{d}\sqrt{2N}\leq (2576+357\log(\varepsilon^{-1}))\sqrt{d}\sqrt{N}
\end{equation*}
which concludes the proof.
By (\ref{157}) which also holds in the case $d=1$ and $N\geq 3$ it is easy to see that for $d=1$ and $m\geq 1$ we observe
\begin{equation}
 \label{256}
 \mathbb{P}(\mathcal{F}_{m,d,\varepsilon})=0.
\end{equation}
Therefore we may assume $d\geq 2$. We now claim
\begin{equation}
 \label{169}
 \mathbb{P}\left(\bigcup_{d\geq 2}\bigcup_{m\geq 1}\mathcal{F}_{m,d.\varepsilon}\right)=\mathbb{P}\left(\bigcup_{d\geq 2}\bigcup_{m\in\{1,\ldots,12\cdot 2^d-1\}}\mathcal{F}_{m,d.\varepsilon}\right).
\end{equation}
Let $d\geq 2$ be given and assume $m\geq 12\cdot 2^d$. Furthermore set $\mu=12\cdot 2^d$. 
By Lemma \ref{168} for $M\in\{2^{m}+1,\ldots,2^{m+1}\}$ we obtain
\begin{equation}
 \label{171}
 MD_M^d(x_{n,i})\leq 2^{\mu}D^d_{2^{\mu}}(x_{n,i})+(M-2^{\mu})D_{2^{\mu},M}^d(x_{n,i}).
\end{equation}
Now observe that since $(x_{2^{\mu},1},\ldots,x_{2^{\mu},d})$ is uniformly distributed the points
\begin{equation*}
 \{(x_{2^{\mu}+1,1},\ldots,x_{2^{\mu}+1,d}),\ldots,(x_{M,1},\ldots,x_{M,d})\}
\end{equation*}
are elements of a randomized Halton sequence denoted by $(q_n)_{n\geq 1}$. Therefore by Lemma \ref{170} and another application of Lemma \ref{168} we have
\begin{equation}
 \label{172}
 \begin{aligned}
  D^d_{2^{\mu},M}(x_{n,i})=D^d_{2^{\mu},M}(q_{n,i}) \:\:\: & \leq \:\:\: \frac{2^{\mu}D^d_{2^{\mu}}(q_{n,i})+MD^d_M(q_{n,i})}{M-2^{\mu}}\\
  & \leq \:\:\: \frac{2^{\mu}\sqrt{d}/\sqrt{2^{\mu}}+M\sqrt{d}/\sqrt{M}}{M-2^{\mu}}\\
  & < \:\:\: \frac{2\sqrt{d}\sqrt{2^{m+1}}}{M-2^{\mu}}.
 \end{aligned}
\end{equation}
Together with (\ref{171}) we get
\begin{equation*}
 \begin{aligned}
  \mathcal{F}_{m,d,\varepsilon}\backslash\mathcal{F}_{\mu-1,d,\varepsilon} \subset & \left\{\max_{M\in\{2^m+1,\ldots,2^{m+1}\}}MD^d_M(x_{n,i})\geq (1821+252\log(\varepsilon^{-1}))\sqrt{d}\sqrt{2^{m+1}}\right\}\mathlarger{\mathlarger{\mathlarger{\mathlarger{\backslash}}}}\\
  & \left\{2^{\mu}D^d_{2^{\mu}}(x_{n,i})\geq (1819+252\log(\varepsilon^{-1}))\sqrt{d}\sqrt{2^{\mu}}\right\}\\
  \subset & \left\{\max_{M\in\{2^m+1,\ldots,2^{m+1}\}}(M-2^{\mu})D^d_{2^{\mu},M}(x_{n,i})> 2\sqrt{d}\sqrt{2^{m+1}}\right\}\\
  = & \emptyset.
 \end{aligned}
\end{equation*}
Therefore for $m\geq \mu$ we have $\mathbb{P}(\mathcal{F}_{m,d,\varepsilon}\backslash\mathcal{F}_{\mu-1,d,\varepsilon})=0$ and (\ref{169}) follows immediately. Thus we may assume $d\geq 2$ and $m\in\{1,\ldots,12\cdot 2^d-1\}$ now.
Furthermore we may assume
\begin{equation}
 \label{179}
 \frac{\sqrt{d}}{\sqrt{2^{m+1}}}\leq\frac{1}{64}  
\end{equation}
since otherwise $\mathcal{F}_{m,d,\varepsilon}=\emptyset$.
Let $\tilde{k}(m)=\max\{k\geq 1:12\cdot 2^k\leq m\}$ and for $m\geq 48$ set $L_m=2^{12\cdot 2^{\tilde{k}(m)}}$ resp. for $m<48$ set $L_m=0$. Moreover we define the sets
\begin{eqnarray*}
 \mathcal{G}_{m,d,\varepsilon} & = &
 \begin{cases}
  \left\{L_mD^d_{L_m}(x_{n,i})\geq (910+126\log(\varepsilon^{-1}))\sqrt{d}\sqrt{L_m}\right\}, & \textnormal{if }L_m>0,\\
  \emptyset, & \textnormal{if }L_m=0,
 \end{cases}\\
 \mathcal{H}_{m,d,\varepsilon} & = & \left\{\max_{L_m+1\leq M\leq 2^{m+1}}(M-L_m)D^d_{L_m,M}(x_{n,i})\geq (909+126\log(\varepsilon^{-1}))\sqrt{d}\sqrt{2^{m+1}}\right\}.
\end{eqnarray*}
Now we claim
\begin{equation}
 \label{173}
 \mathcal{F}_{m,d,\varepsilon}\subseteq\mathcal{G}_{m,d,\varepsilon}\cup\mathcal{H}_{m,d,\varepsilon}
\end{equation}
for all $d\geq 2$ and $m\in\{1,\ldots,12\cdot 2^d-1\}$.
Since this trivially holds for $L_m=0$ we may assume $L_m\geq 1$ and therefore we have $m\geq 48$ and $\tilde{k}(m)\geq 2$. By Lemma \ref{168} for the complement of $\mathcal{G}_{m,d,\varepsilon}\cup\mathcal{H}_{m,d,\varepsilon}$ we observe
\begin{eqnarray*}
 \max_{L_m+1\leq M\leq 2^{m+1}}MD^d_M(x_{n,i}) & \leq & \max_{L_m+1\leq M\leq 2^{m+1}}(L_mD^d_{L_m}(x_{n,i})+(M-L_m)D^d_{{L_m},M}(x_{n,i}))\\
 & \leq & (910+126\log(\varepsilon^{-1}))\sqrt{d}\sqrt{L_m}\\
 && +(909+126\log(\varepsilon^{-1}))\sqrt{d}\sqrt{2^{m+1}}\\
 & \leq & (1819+252\log(\varepsilon^{-1}))\sqrt{d}\sqrt{2^{m+1}}.
\end{eqnarray*}
Thus we have (\ref{173}). Now for any $d\geq 2$ we shall show
\begin{equation}
 \label{174}
 \bigcup_{m\in\{1,\ldots,12\cdot 2^d-1\}}\mathcal{F}_{m,d,\varepsilon}\subseteq\bigcup_{m\in\{1,\ldots,12\cdot 2^d-1\}}\mathcal{H}_{m,d,\varepsilon}.
\end{equation}
For any $k\geq 2$ by definition of $L_m$ we have
\begin{equation*}
 \mathcal{G}_{12\cdot 2^k,m,\varepsilon}=\mathcal{G}_{12\cdot 2^k+1,m,\varepsilon}=\cdots=\mathcal{G}_{12\cdot 2^{k+1}-1,m,\varepsilon}.
\end{equation*}
Therefore by (\ref{173}) we obtain
\begin{equation}
 \label{175}
 \bigcup_{m\in\{1,\ldots,12\cdot 2^d-1\}}\mathcal{F}_{m,d,\varepsilon}\subseteq\bigcup_{k\in\{2,\ldots,d-1\}}\mathcal{G}_{12\cdot 2^k,d,\varepsilon}\cup\bigcup_{m\in\{1,\ldots,12\cdot 2^d\}}\mathcal{H}_{m,d,\varepsilon}.
\end{equation}
For $k=2$ we have $m=12\cdot 2^2=48$ and $L_m=2^{48}$. With $L_{47}=0$ we get
\begin{eqnarray*}
 \mathcal{G}_{48,d,\varepsilon} & \subseteq & \left\{2^{48}D^d_{2^{48}}(x_{n,i})\geq (910+126\log(\varepsilon^{-1}))\sqrt{d}\sqrt{2^{48}}\right\}\\
 & \subseteq & \left\{L_{47}+(2^{48}-L_{47})D^d_{L_{47},2^{48}}(x_{n,i})\geq (910+126\log(\varepsilon^{-1}))\sqrt{d}\sqrt{2^{48}}\right\}\\
 & \subseteq & \left\{(2^{48}-L_{47})D^d_{L_{47},2^{48}}(x_{n,i})\geq (909+126\log(\varepsilon^{-1}))\sqrt{d}\sqrt{2^{48}}\right\}\\
 & \subseteq & \mathcal{H}_{47,d,\varepsilon}
\end{eqnarray*}
where the second line follows by Lemma \ref{168}.
For $k\geq 3$ and $m=12\cdot 2^k$ we have $\tilde{k}(m)=k$ and $L_m=2^{12\cdot 2^k}$. Moreover we obtain $\tilde{k}(m-1)=k-1$ and furthermore $L_{m-1}=2^{12\cdot 2^{k-1}}=\sqrt{L_m}$.
Thus we have
\begin{eqnarray*}
 \mathcal{G}_{m,d,\varepsilon} & \subseteq & \left\{L_mD^d_{L_m}(x_{n,i})\geq (910+126\log(\varepsilon^{-1}))\sqrt{d}\sqrt{L_m}\right\}\\
 & \subseteq & \left\{L_{m-1}+(L_m-L_{m-1})D^d_{L_{m-1},L_m}(x_{n,i})\geq (910+126\log(\varepsilon^{-1}))\sqrt{d}\sqrt{L_m}\right\}\\
 & \subseteq & \left\{(L_m-L_{m-1})D^d_{L_{m-1},L_m}(x_{n,i})\geq (909+126\log(\varepsilon^{-1}))\sqrt{d}\sqrt{L_m}\right\}\\
 & \subseteq & \mathcal{H}_{m-1,d,\varepsilon}.
\end{eqnarray*}
Together with (\ref{175}) we observe (\ref{174}). Thus by (\ref{167}), (\ref{256}) and (\ref{169}) the Theorem is proved if we show
\begin{equation}
 \label{178}
 \sum_{d\geq 2}\sum_{m\in\{1,\ldots,12\cdot 2^d-1\}}\mathbb{P}(\mathcal{H}_{m,d,\varepsilon})\leq\varepsilon.
\end{equation}
Now we shall prove
\begin{equation}
 \label{202}
 \mathbb{P}(\mathcal{H}_{m,d,\varepsilon})\leq \frac{\varepsilon}{6\cdot 2^{2d}}
\end{equation}
for all $d\geq 2$ and $m\in\{1,\ldots,12\cdot 2^d-1\}$. Then (\ref{178}) follows by
\begin{equation*}
 \sum_{d\geq 2}\sum_{m\in\{1,\ldots,12\cdot 2^d-1\}}\mathbb{P}(\mathcal{H}_{m,d,\varepsilon})\leq \sum_{d\geq 2}12\cdot 2^d\cdot\frac{\varepsilon}{6}\cdot 2^{-2d}=\varepsilon.
\end{equation*}
To prove (\ref{202}) let $d\geq 2$ and $m\in\{1,\ldots,12\cdot 2^d-1\}$ be fixed now. To estimate $D^d_{L_m,M}(x_{n,i})$ we define a finite system of subsets of $[0,1)^d$ with the help of $\delta$-bracketing covers such that $[0,y)$ for any $y\in [0,1)^d$ can be
approximated well enough by a union of this sets. Set
\begin{equation}
 \label{180}
 H=\left\lceil\frac{m+1}{2}-\frac{\log_2(d)}{2}-2\right\rceil.
\end{equation}
As a consequence for any $h\in\{0,\ldots,H\}$ we have
\begin{equation}
 \label{181}
 \sqrt{d}\sqrt{2^{m+1}}\leq 2^{m-h}.
\end{equation}
For any $h\in\{1,\ldots,H\}$ let $\Delta_h$ be a $2^{-h}$-bracketing cover of $[0,1)^d$. By Corollary \ref{177} we may assume
\begin{equation}
 |\Delta_h|\leq \frac{1}{2}(2e)^d(2^{h+2}+1)^d.
\end{equation}
For any $y\in [0,1)^d$ we now define a finite sequence of points $\beta_h(y)$ for $h\in\{0,\ldots,H+1\}$ in the following manner. Let $(v,w)\in\Delta_H$ be such that $v\leq y\leq w$. We set $\beta_{H+1}(y)=w$ and $\beta_H(y)=v$.
The points $\beta_1(y),\ldots,\beta_{H-1}(y)$ are defined by induction. Thus assume that for some $h\in\{1,\ldots,H-1\}$ the point $\beta_{h+1}(y)$ is already defined. Let $(v,w)\in\Delta_{h}$ with $v\leq \beta_{h+1}(y)\leq w$ and set $\beta_{h}(y)=v$.
Moreover set $\beta_0(y)=0$. Therefore we observe
\begin{equation*}
 0=\beta_0(y)\leq \beta_1(y)\leq \cdots \leq \beta_H(y)\leq x\leq \beta_{H+1}(y)\leq 1.
\end{equation*}
For $h\in\{0,\ldots,H-1\}$ we have $(\beta_h(y),w)\in\Delta_h$ for some point $w\in [0,1)^d$.
Furthermore we have $(\beta_H(y),\beta_{H+1}(y))\in\Delta_H$.
Then by Corollary \ref{177} for $h\in\{0,\ldots,H+1\}$ and $i\in\{1,\ldots, d\}$ there exist integers $a_{h,i}\in\{0,\ldots,2^{(\lceil \log_2(i)+1\rceil)(h+1)}\}$ such that
\begin{equation}
 \label{182}
 (\beta_h(y))_i=2^{-(\lceil \log_2(i)+1\rceil)(h+1)}a_{h,i}.
\end{equation}
For $h\in\{0,\ldots,H\}$ set $K_h(y)=\overline{[\beta_h(y),\beta_{h+1}(y))}$. Note that the sets $K_h(y)$ are pairwise disjoint and satisfy
\begin{equation}
 \label{183}
 \bigcup_{h=0}^{H-1}K_h(y)\subseteq [0,x)\subseteq\bigcup_{h=0}^{H}K_h(y)
\end{equation}
By definition $\beta_h(y)\leq \beta_{h+1}(y)\leq w$ for some $w\in [0,1)^d$ with $(\beta_h(y),w)\in\Delta_h$ and hence
\begin{equation}
 \label{184}
 \lambda(K_h(y))\leq \lambda\left(\overline{[\beta_h(y),w)}\right)\leq 2^{-h}
\end{equation}
for any $h\in\{0,\ldots,H\}$.
Now define 
\begin{equation*}
S_h=\left\{\overline{[\beta_h(y),\beta_{h+1}(y))}:x\in [0,1)^d\right\}. 
\end{equation*}
Observe that we may define the points $\beta_h$ such that $\beta_h(y)=\beta_h(z)$ for $y,z\in [0,1)^d$ with $\beta_{h+1}(y)=\beta_{h+1}(z)$. Therefore by Corollary \ref{177} we have
\begin{equation}
 \label{185}
 |S_h|=\left|\left\{\beta_{h+1}(y):y\in [0,1)^d\right\}\right|\leq |\Delta_{h+1}|\leq \frac{1}{2}(2e)^d(\sqrt{5})^{(h+3)d}
\end{equation}
for any integer $h\in\{0,\ldots,H\}$.
For $m\geq 48$ we set $s=\tilde{k}(m)$. Otherwise we set $s=1$. Let now $n\in\{L_m+1,\ldots,2^{m+1}\}$ be an integer. For $m<48$ we have $s=1$ and therefore we obtain $x_{n,i}=T_{p_i}(x_{n-1,i})$ for $i\leq s$ by definition while for $i>s$ we get
\begin{equation*}
 n\leq 2^{m+1}\leq 2^{48}\leq 2^{12\cdot 2^i}.
\end{equation*}
Thus for $i\geq 2$ we have $x_{n,i}=\langle 2^{\lceil\log_2(i)\rceil+1}x_{n-1,i}\rangle$. For $m\geq 48$ and $i\leq s=\tilde{k}(m)$ we get
\begin{equation*}
 n>L_m=2^{12\cdot 2^{\tilde{k}(m)}}\geq 2^{12\cdot 2^i}
\end{equation*}
and thus we obtain $x_{n,i}=T_{p_i}(x_{n-1,i})$. Furthermore for $i>s=\tilde{k}(m)$ we observe
\begin{equation*}
 n\leq 2^{m+1}\leq 2^{12\cdot 2^i}
\end{equation*}
and we obtain $x_{n,i}=\langle 2^{\lceil\log_2(i)\rceil+1}x_{n-1,i}\rangle$.
We see that in the sequence
\begin{equation*}
 \left\{(x_{L_m+1,1},\ldots,x_{L_m+1,d}),\ldots,(x_{2^{m+1},1},\ldots,x_{2^{m+1},d})\right\}
\end{equation*}
the first $s$ coordinates form a randomized Halton sequence while the sequence formed by the remaining coordinates is a sequence of fractional parts of the product of some initial value and elements of a lacunary sequence.
Hence for any $M\in\{2^{m}+1,\ldots,2^{m+1}\}$ by Lemma \ref{168} and \ref{170} we have
\begin{equation}
 \label{188}
 \begin{aligned}
  (M-L_m)D^s_{L_m,M}(x_{n,i}) \:\:\: & \leq \:\:\: L_mD^s_{L_m}(x_{n,i})+MD^s_M(x_{n,i})\\
   & \leq \:\:\: \sqrt{s}\sqrt{L_m}+\sqrt{s}\sqrt{M}\leq 2\sqrt{s}\sqrt{M}.
 \end{aligned}
\end{equation}
For some $h\in\{1,\ldots,H+1\}$ and a point $y\in [0,1)^d$ the point $\beta_h(y)$ can be written as $(u_h(y),v_h(y))$ for $u_h(y)\in [0,1)^s$ and $v_h(y)\in [0,1)^{d-s}$. Moreover set $U_h(y)=[0,u_h(y))$ and $V_h(y)=[0,v_h(y))$. Thus we have $U_h(y)\times V_h(y)=[0,\beta_h(y))$.
Observe that any set $K_h(y)\in S_h$ may be written as
\begin{equation*}
 \begin{aligned}
  K_h(y) \:\:\: & = \:\:\: \,\, \overline{[\beta_h(y),\beta_{h+1}(y))}\\
  & = \:\:\: ((U_{h+1}(y)\backslash U_h(y))\times V_{h+1}(y))\cup (U_h(y)\times(V_{h+1}(y)\backslash V_h(y))).
 \end{aligned}
\end{equation*}
For $h=0$ we simply have $K_0(y)=U_1(y)\times V_1(y)$. Furthermore set $U_0(y)=V_0(y)=\emptyset$. Thus by (\ref{184}) we observe
\begin{equation}
 \label{192}
 \lambda(U_{h+1}(y)\backslash U_h(y))\cdot\lambda (V_{h+1}(y))+\lambda(U_h(y))\cdot\lambda (V_{h+1}(y)\backslash V_h(y)) \leq \lambda(K_h(y))\leq 2^{-h}
\end{equation}
for $h\in\{0,\ldots,H\}$.
Now let $y\in [0,1)^d$ be some arbitrary fixed point. Note that hereafter we skip the point $y$ in the notation of the points $\beta_h$ and the sets $K_h$ resp. $U_h$ and $V_h$ to simplify notations. Furthermore let $L_m+1\leq M\leq 2^{m+1}$ be an integer. For simplicity we write $q_n=(x_{n,1},\ldots,x_{n,s})$ and $r_n=(x_{n,s+1},\ldots,x_{n,d})$. Then by (\ref{183}) we have
\begin{equation}
 \label{186}
 \begin{aligned}
  \sum_{n=L_m+1}^M\mathbf{1}_{[0,y)}(x_n) \geq & \sum_{n=L_m+1}^{M}\mathbf{1}_{[0,\beta_H)}(x_n)\\
  = & \sum_{n=L_m+1}^M\mathbf{1}_{U_H}(q_n)\cdot \mathbf{1}_{V_H}(r_n)\\
  = & \sum_{n=L_m+1}^M\mathbf{1}_{U_1}(q_n)\cdot \mathbf{1}_{V_1}(r_n)\\
  & +\sum_{h=1}^{H-1}\sum_{n=L_m+1}^M\left(\mathbf{1}_{U_{h+1}\backslash U_h}(q_n)\cdot \mathbf{1}_{V_{h+1}}(r_n)+\mathbf{1}_{U_h}(q_n)\cdot \mathbf{1}_{V_{h+1}\backslash V_h}(r_n)\right).
 \end{aligned}
\end{equation}
Analogously we also get
\begin{equation}
 \label{187}
 \begin{aligned}
  \sum_{n=L_m+1}^M\mathbf{1}_{[0,y)}(x_n)\leq & \sum_{n=L_m+1}^M\mathbf{1}_{U_1}(q_n)\cdot \mathbf{1}_{V_1}(r_n)\\
  & +\sum_{h=1}^{H}\sum_{n=L_m+1}^M\left(\mathbf{1}_{U_{h+1}\backslash U_h}(q_n)\cdot \mathbf{1}_{V_{h+1}}(r_n)+\mathbf{1}_{U_h}(q_n)\cdot \mathbf{1}_{V_{h+1}\backslash V_h}(r_n)\right).
 \end{aligned}
\end{equation}
By using maximal Bernstein inequality we now shall give a lower bound on the probability that the system of inequalities
\begin{eqnarray}
 \label{204}
  \max_{L_m+1\leq M\leq 2^{m+1}}\left|\sum_{n=L_m+1}^M\mathbf{1}_{U_{h+1}\backslash U_h}(q_n)\mathbf{1}_{V_{h+1}}(r_n)-\mathbf{1}_{U_{h+1}\backslash U_h}(q_n)\lambda(V_{h+1})\right| & > & t,\quad\quad\quad\quad\\
 \label{205}
  \max_{L_m+1\leq M\leq 2^{m+1}}\left|\sum_{n=L_m+1}^M\mathbf{1}_{U_h}(q_n)\mathbf{1}_{V_{h+1}\backslash V_h}(r_n)-\mathbf{1}_{U_h}(q_n)\lambda(V_{h+1}\backslash V_h)\right| & > & t,\\
 \label{206}
  \max_{L_m+1\leq M\leq 2^{m+1}}\left|\sum_{n=L_m+1}^M\mathbf{1}_{U_1}(q_n)\mathbf{1}_{V_1}(r_n)-\mathbf{1}_{U_1}(q_n)\lambda(V_1)\right| & > & t
\end{eqnarray}
holds for all sets $U_h$, $U_{h+1}$, $V_h$ and $V_{h+1}$ with $h\in\{1,\ldots,H\}$ and some $t>0$ to specified later. Set $\kappa=\kappa_h=\lceil \log_2(h+2)\rceil$.
By Lemma \ref{170} and (\ref{188}) for any $h\in\{1,\ldots,H\}$ we have
\begin{equation}
 \label{189}
 \begin{aligned}
  \sum_{n=L_m+1}^{M}\mathbf{1}_{U_{h+1}\backslash U_h}(q_n) \:\:\: & \leq \:\:\: \sqrt{s}\cdot\sqrt{\frac{M-L_m}{2^{\kappa}}+1}+\sum_{\substack{n\in\{L_m+1,\ldots,2^{m+1}\},\\n\equiv \gamma \pmod{2^{\kappa}}}}\lambda(U_{h+1}\backslash U_h)\\
  & \leq \:\:\: \left(2^{m+1-\kappa}+1\right)\lambda(U_{h+1}\backslash U_h)+\sqrt{s}\cdot\sqrt{2^{m+1-\kappa}+1}
 \end{aligned}
\end{equation}
and
\begin{equation}
 \label{190}
 \sum_{n=L_m+1}^{M}\mathbf{1}_{U_h}(q_n)\leq \left(2^{m+1-\kappa}+1\right)\lambda(U_{h+1}\backslash U_h)+\sqrt{s}\cdot\sqrt{2^{m+1-\kappa}+1}.
\end{equation}
Now let $h\in\{0,\ldots,H\}$ be fixed and set $\mathcal{A}_h=V_{h+1}\backslash V_h$ resp. $\mathcal{A}_h=V_{h+1}$. Furthermore define by  $f_{\mathcal{A}_h}(x)=f_{\mathcal{A}_h}(r_n)=\mathbf{1}_{\mathcal{A}_h}(x)-\lambda(\mathcal{A}_h)$ a real-valued function on $[0,1)^{d-s}$. 
We now shall show that for any system of indices $n_1,\ldots,n_k$ with $n_{l+1}-n_l\geq h+2$ for all $l\in\{1,\ldots,k-1\}$ the random variables $f_{\mathcal{A}_h}(r_n)$ are stochastically independent, i.e.
\begin{equation}
 \label{194}
 \mathbb{P}\left(f_{\mathcal{A}_h}(r_{n_1})=c_1,\ldots,f_{\mathcal{A}_h}(r_{n_k})=c_k\right)=\prod_{l=1}^k\mathbb{P}\left(f_{\mathcal{A}_h}(r_{n_l})=c_l\right).
\end{equation}
We only prove the case $k=2$. The general case follows by induction. By (\ref{182}) the set $\mathcal{A}_h$ is a union of axis-parallel boxes such that each corner of any box is of the form
\begin{equation}
 \label{193}
 \left(2^{-(\lceil\log_2(s+1)\rceil+1)(h+2)}a_{s+1},\ldots,2^{-(\lceil\log_2(d)\rceil+1)(h+2)}a_{d}\right)
\end{equation}
such that $a_i\in\{0,1,\ldots,2^{(\lceil\log_2(d)\rceil+1)(h+2)}\}$ for any $i\in\{s+1,\ldots,d\}$.
Furthermore let $n,n'\in\{L_m+1,\ldots,M\}$ be two indices with $n'-n\geq h+2$. We define a decomposition of $[0,1)^{d-s}$ by
\begin{multline*}
 \Sigma=\left\{\prod_{i=s+1}^d\left[2^{-(\lceil\log_2(i)\rceil+1)n'}a_i,2^{-(\lceil\log_2(i)\rceil+1)n'}(a_i+1)\right):\right.\\
 \left.a_i\in\left\{0,1,\ldots,2^{(\lceil\log_2(i)\rceil+1)n'}-1\right\},i\in\{s+1,\ldots,d\}\right\}.
\end{multline*}
Note that by (\ref{193}) the function $f_{\mathcal{A}_h}$ is constant on any box $\mathcal{B}\in\Sigma$. For some $c_1\in\mathbb{R}$ define
\begin{equation*}
 \Sigma_{c_1}=\left\{\mathcal{B}\in\Sigma:f_{\mathcal{A}_h}(r_{n})=c_1 \textnormal{ for all }r_1=(r_{1,s+1},\ldots,r_{1,d})\in\mathcal{B}\right\}.
\end{equation*}
Since $x_{n',i}=2^{(\lceil\log_2(i)\rceil+1)(n'-1)}x_{1,i}$ for all $i\in\{s+1,\ldots,d\}$ we have $f_{\mathcal{A}_h}(r_{n'})=f_{\mathcal{A}_h}(r'_{n'})$ where $r'_{n'}=(x'_{n',s+1},\ldots,x'_{n',d})$ with $x'_{n,i}=2^{(\lceil\log_2(i)\rceil+1)(n'-1)}x'_{1,i}$ is an
instance of the matrix for some initial value $r'_1=(x'_{1,s+1},\ldots,x'_{1,d})$ with $x'_{1,i}=x_{1,i}+2^{-(\lceil\log_2(i)\rceil+1)(n'-1)}a_i$ and $a_i\in\{0,1,\ldots,2^{(\lceil\log_2(i)\rceil+1)(n'-1)}-1\}$ for all $i\in\{s+1,\ldots,d\}$.
Therefore for any $c_2\in\mathbb{R}$ and any $\mathcal{B},\mathcal{B}'\in\Sigma$ we have
\begin{equation*}
 \mathbb{P}\left(f_{\mathcal{A}_h}(r_{n'})=c_2|r_1\in\mathcal{B}\right)=\mathbb{P}\left(f_{\mathcal{A}_h}(r_{n'})=c_2|r_1\in\mathcal{B}'\right).
\end{equation*}
Hence for any $c_2\in\mathbb{R}$ and any $\mathcal{B}\in\Sigma$ we get
\begin{eqnarray*}
 \mathbb{P}\left(f_{\mathcal{A}_h}(r_{n'})=c_2\right) & = & \sum_{\mathcal{B}'\in\Sigma}\mathbb{P}\left(f_{\mathcal{A}_h}(r_{n'})=c_2|r_1\in\mathcal{B}'\right)\mathbb{P}(r_1\in\mathcal{B}')\\
 & = & \mathbb{P}\left(f_{\mathcal{A}_h}(r_{n'})=c_2|r_1\in\mathcal{B}\right)\sum_{\mathcal{B}'\in\Sigma}\mathbb{P}(r_1\in\mathcal{B}')\\
 & = & \mathbb{P}\left(f_{\mathcal{A}_h}(r_{n'})=c_2|r_1\in\mathcal{B}\right).
\end{eqnarray*}
Moreover for any $c_1,c_2\in\mathbb{R}$ we obtain
\begin{multline*}
 \mathbb{P}\left(f_{\mathcal{A}_h}(r_{n'})=c_2|f_{\mathcal{A}_h}(r_n)=c_1\right)\\
 \begin{aligned}
  = & \frac{\mathbb{P}\left(f_{\mathcal{A}_h}(r_{n'})=c_2,f_{\mathcal{A}_h}(r_n)=c_1\right)}{\mathbb{P}\left(f_{\mathcal{A}_h}(r_n)=c_1\right)}\\
  = & \frac{\sum_{\mathcal{B}\in\Sigma}\mathbb{P}\left(f_{\mathcal{A}_h}(r_{n'})=c_2,f_{\mathcal{A}_h}(r_n)=c_1|r_1\in\mathcal{B}\right)\mathbb{P}(r_1\in\mathcal{B})}{\mathbb{P}\left(f_{\mathcal{A}_h}(r_n)=c_1\right)}\\
  = & \sum_{\mathcal{B}\in\Sigma_{c_1}}\mathbb{P}\left(f_{\mathcal{A}_h}(r_{n'})=c_2|r_1\in\mathcal{B}\right)\frac{\mathbb{P}(r_1\in\mathcal{B})}{\mathbb{P}\left(f_{\mathcal{A}_h}(r_n)=c_1\right)}\\
  = & \mathbb{P}\left(f_{\mathcal{A}_h}(r_{n'})=c_2\right).
 \end{aligned}
\end{multline*}
Thus (\ref{194}) is proved.
Furthermore set 
\begin{eqnarray*}
 Q(L_m,M,\gamma) & = & \left\{n\in\left\{L_m+1,\ldots, M\right\}:q_n\in U_{h+1}\backslash U_h,n\equiv \gamma \,\, \pmod {2^{\kappa}}\right\},\\
 Q'(L_m,M,\gamma) & = & \left\{n\in\left\{L_m+1,\ldots, M\right\}:q_n\in U_h,n\equiv \gamma \,\, \pmod {2^{\kappa}}\right\}.
\end{eqnarray*}
Then for $h\in \{1,\ldots, H\}$ by Lemma \ref{191} we have
\begin{multline*}
 \mathbb{P}\left(\max_{M\in\{L_m+1,\ldots,2^{m+1}\}}\left|\sum_{n=L_m+1}^M\mathbf{1}_{U_{h+1}\backslash U_h}(q_n)\cdot\mathbf{1}_{V_{h+1}}(r_n)-\mathbf{1}_{U_{h+1}\backslash U_h}(q_n)\lambda(V_{h+1})\right|>t\right)\\
 \begin{aligned}
  \leq & \sum_{\gamma=1}^{2^{\kappa}}\mathbb{P}\left(\max_{n\in\left\{L_m+1,\ldots, M\right\}}\left|\sum_{n\in Q(L_m,M,\gamma)}\mathbf{1}_{V_{h+1}}(r_n)-\lambda(V_{h+1})\right|>\frac{t}{2^{\kappa}}\right)\\
  \leq & 2\sum_{\gamma=1}^{2^{\kappa}}\exp\left(-\frac{t^2/2^{2\kappa}}{2\left(\sum_{n\in Q(L_m,M,\gamma)}1\right)\lambda(V_{h+1})(1-\lambda(V_{h+1}))+2t/(3\cdot 2^{\kappa})}\right).
 \end{aligned}
\end{multline*}
Thus by (\ref{189}) we obtain
\begin{multline*}
 \mathbb{P}\left(\max_{M\in\{L_m+1,\ldots,2^{m+1}\}}\left|\sum_{n=L_m+1}^M\mathbf{1}_{U_{h+1}\backslash U_h}(q_n)\cdot\mathbf{1}_{V_{h+1}}(r_n)-\mathbf{1}_{U_{h+1}\backslash U_h}(q_n)\lambda(V_{h+1})\right|>t\right)\\
 \leq 2^{\kappa+1}\exp\left(-\frac{t^2/2^{1.5\kappa}}{2^{m+3}\lambda(U_{h+1}\backslash U_h)\lambda(V_{h+1})+2\sqrt{2}\cdot\sqrt{s}\cdot\sqrt{2^{m+1}}\lambda(V_{h+1})+2t/3}\right).
\end{multline*}
Furthermore (\ref{181}) and (\ref{192}) yield
\begin{multline}
 \label{197}
 \mathbb{P}\left(\max_{M\in\{L_m+1,\ldots,2^{m+1}\}}\left|\sum_{n=L_m+1}^M\mathbf{1}_{U_{h+1}\backslash U_h}(q_n)\cdot\mathbf{1}_{V_{h+1}}(r_n)-\mathbf{1}_{U_{h+1}\backslash U_h}(q_n)\lambda(V_{h+1})\right|>t\right)\\
 \leq 2^{\kappa+1}\exp\left(-\frac{t^2/2^{1.5\kappa}}{(8+2\sqrt{2})\cdot 2^{m-h}+2t/3}\right).
\end{multline}
Similarly using (\ref{190}) we get
\begin{multline}
 \label{198}
 \mathbb{P}\left(\max_{M\in\{L_m+1,\ldots,2^{m+1}\}}\left|\sum_{n=L_m+1}^M\mathbf{1}_{U_h}(q_n)\cdot\mathbf{1}_{V_{h+1}\backslash V_h}(r_n)-\mathbf{1}_{U_h}(q_n)\lambda(V_{h+1}\backslash V_h)\right|>t\right)\\
 \begin{aligned}
  \leq & 2^{\kappa+1}\exp\left(-\frac{t^2/2^{1.5\kappa}}{(8+2\sqrt{2})\cdot 2^{m-h}+2t/3}\right).
 \end{aligned}
\end{multline}
Now set $t=C_1\sqrt{d}\sqrt{2^{m+1}}\sqrt{h\cdot 2^{1.5\kappa-h}}$ for a constant $C_1>0$ to specified later. Observe that by (\ref{181}) we have $t\leq 2^{m-h+1}C_1$.

Therefore by (\ref{197}) we get
\begin{multline}
 \label{199}
 \mathbb{P}\left(\max_{M\in\{L_m+1,\ldots,2^{m+1}\}}\left|\sum_{n=L_m+1}^M\mathbf{1}_{U_{h+1}\backslash U_h}(q_n)\cdot\mathbf{1}_{V_{h+1}}(r_n)-\mathbf{1}_{U_{h+1}\backslash U_h}(q_n)\lambda(V_{h+1})\right|>t\right)\\
 \begin{aligned}
  \leq & 4(h+2)\exp\left(-\frac{\left(C_1\sqrt{d}\sqrt{2^{m+1}}\sqrt{h\cdot 2^{1.5\kappa-h}}\right)^2}{2^{1.5\kappa}(8+2\sqrt{2}+2C_1)2^{m-h}}\right)\\
  \leq & 4\exp\left(-\left(\frac{2C_1^2}{8+2\sqrt{2}+2C_1}-1\right)hd\right)
 \end{aligned}
\end{multline}
where the last line follows by $(h+2)\leq e^{hd}$ for $d\geq 2$. Similarly using (\ref{198}) we have
\begin{multline}
 \label{200}
 \mathbb{P}\left(\max_{M\in\{L_m+1,\ldots,2^{m+1}\}}\left|\sum_{n=L_m+1}^M\mathbf{1}_{U_h}(q_n)\cdot\mathbf{1}_{V_{h+1}\backslash V_h}(r_n)-\mathbf{1}_{U_h}(q_n)\lambda(V_{h+1}\backslash V_h)\right|>t\right)\\
  \leq 4\exp\left(-\left(\frac{2C_1^2}{8+2\sqrt{2}+2C_1}-1\right)hd\right).
\end{multline}
For $h=0$ set $t=C_2\sqrt{d}\sqrt{2^{m+1}}$ for some constant $C_2>0$ to be specified later. Thus by using a similar argumentation as above we get
\begin{multline}
 \label{201}
 \mathbb{P}\left(\max_{M\in\{L_m+1,\ldots,2^{m+1}\}}\left|\sum_{n=L_m+1}^M\mathbf{1}_{U_1}(q_n)\mathbf{1}_{V_1}(r_n)-\mathbf{1}_{U_1}(q_n)\lambda(V_1)\right|>t\right)\\
  \leq  4\exp\left(-\frac{t^2/2}{(8+2\sqrt{2})\cdot 2^{m}+2t/3}\right)         \leq  4\exp\left(-\frac{C_2^2}{8+2\sqrt{2}+2/3\cdot C_2}d\right).
\end{multline}
Define
\begin{equation}
 \label{203}
 C_3=\frac{2C_1^2}{8+2\sqrt{2}+2C_1}-1, \quad C_4=\frac{C_2^2}{8+2\sqrt{2}+2/3\cdot C_2}.
\end{equation}
Observe that by (\ref{185}) and sufficiently large constants $C_1$, $C_2$ resp. $C_3$, $C_4$ the system of inequalities (\ref{204}), (\ref{205}) and (\ref{206}) hold on a set of measure which is bounded from below by
\begin{equation}
 \label{207}
 1-\frac{1}{2}(2e)^d(\sqrt{5})^{3d}\cdot 4e^{-C_4d}-(2e)^d\sum_{h=1}^H(\sqrt{5})^{(h+3)d}\cdot 4e^{-C_3d}\geq 1-\frac{\varepsilon}{6\cdot 2^{2d}}.
\end{equation}
Now we shall find some constants $C_1$ and $C_2$ such that (\ref{207}) is true. It is easy to see that for $C_3\geq 2.7$ we have
\begin{equation*}
 4(50e)^d\sum_{h=1}^H(\sqrt{5})^{(h-1)d}e^{-C_3hd}=4(50e)^de^{-C_3d}\sum_{h=1}^H(\sqrt{5}e^{-C_3})^{(h-1)d}\leq 4.1(50e)^de^{-C_3d}.
\end{equation*}
Therefore we can estimate the left-hand side of (\ref{207}) by
\begin{multline*}
 1-\frac{1}{2}(2e)^d(\sqrt{5})^{3d}\cdot 4e^{-C_4d}-(2e)^d\sum_{h=1}^H(\sqrt{5})^{(h+3)d}\cdot 4e^{-C_3d}\\
 \begin{aligned}
  \geq & 1-\left(2/(\sqrt{5})^d+4.1\right)(50e)^de^{-\min(C_3,C_4)d}\\
  \geq & 1-4.5\cdot e^{(1+\log(50)-\min(C_3,C_4))d}.
 \end{aligned}
\end{multline*}
Thus (\ref{207}) holds if
\begin{equation*}
 -\log(4.5)+\left(\min(C_3,C_4)-1-\log(50)\right)d\geq \log(\varepsilon^{-1})+\log(6)+2d.
\end{equation*}
By $d\geq 2$ it can easily be shown that (\ref{207}) is true for
\begin{equation}
 \label{208}
 \min(C_3,C_4)\geq 7.947+\frac{\log(\varepsilon^{-1})}{2}.
\end{equation}
By (\ref{203}) this holds for 
\begin{equation*}
 C_2\geq \frac{15.894+\log(\varepsilon^{-1})}{6}+\sqrt{86.054+(4+\sqrt{2})\log(\varepsilon^{-1})+\left(\frac{15.894+\log(\varepsilon^{-1})}{6}\right)^2}
\end{equation*}
and because of $\sqrt{A+B}\leq\sqrt{A}+\sqrt{B}$ we may choose
\begin{equation}
 \label{209}
 C_2=14.575+5.748\log(\varepsilon^{-1}).
\end{equation}
Similarly (\ref{208}) holds for
\begin{equation*}
 C_1\geq \frac{17.894+\log(\varepsilon^{-1})}{4}+\sqrt{48.441+2.708\log(\varepsilon^{-1})+\left(\frac{17.894+\log(\varepsilon^{-1})}{4}\right)^2}.
\end{equation*}
Thus we may take
\begin{equation}
 \label{210}
 C_1=15.907+2.146\log(\varepsilon^{-1}).
\end{equation}
Therefore for any $M\in\{L_m+1,\ldots,2^{m+1}\}$ by using (\ref{187}), (\ref{199}), (\ref{200}), (\ref{201}), (\ref{203}), (\ref{207}), (\ref{209}) and (\ref{210}) we get
\begin{multline*}
 \sum_{n=L_m+1}^M\mathbf{1}_{[0,y)}(x_n)\\
 \begin{aligned}
  \leq & \sum_{n=L_m+1}^M\mathbf{1}_{U_1}(q_n)\cdot \mathbf{1}_{V_1}(r_n)\\
  & +\sum_{h=1}^{H}\sum_{n=L_m+1}^M\left(\mathbf{1}_{U_{h+1}\backslash U_h}(q_n)\cdot \mathbf{1}_{V_{h+1}}(r_n)+\mathbf{1}_{U_h}(q_n)\cdot \mathbf{1}_{V_{h+1}\backslash V_h}(r_n)\right)\\
  \leq & \sum_{n=L_m+1}^M\lambda(U_{H+1})\lambda(V_{H+1})+\sum_{n=L_m+1}^M\lambda(V_{H+1})\left(\mathbf{1}_{U_{H+1}}(q_n)-\lambda(U_{H+1})\right)\\
  & + 2(15.907+2.146\log(\varepsilon^{-1}))\sqrt{d}\sqrt{2^{m+1}}\sum_{h=1}^H\sqrt{h\cdot 2^{1.5(1+\log_2(h+2))-h}}\\
  & +(14.575+5.748\log(\varepsilon^{-1}))\sqrt{d}\sqrt{2^{m+1}}
 \end{aligned}
\end{multline*}
with probability at least $1-\varepsilon/6\cdot 2^{-2d}$.
Thus with 
\begin{equation*}
 \sum_{h=1}^H\sqrt{h\cdot 2^{1.5(1+\log_2(h+2))-h}}\leq 27.917
\end{equation*} 
we obtain
\begin{multline*}
 \sum_{n=L_m+1}^M\mathbf{1}_{[0,y)}(x_n)\\
 \begin{aligned}
  \leq & \sum_{n=L_m+1}^M\lambda(U_{H+1})\lambda(V_{H+1})+\sum_{n=L_m+1}^M\lambda(V_{H+1})\left(\mathbf{1}_{U_{H+1}}(q_n)-\lambda(U_{H+1})\right)\\
  & +(902.726+125.568\log(\varepsilon^{-1}))\sqrt{d}\sqrt{2^{m+1}}.
 \end{aligned}
\end{multline*}
By (\ref{180}) and (\ref{188}) we have
\begin{multline}
 \label{211}
 \sum_{n=L_m+1}^M\mathbf{1}_{[0,y)}(x_n)\\
 \begin{aligned}
  \leq & (M-L_m)\lambda([0,\beta_{H+1}(x))+2\sqrt{d}\sqrt{2^{m+1}}+(902.726+125.568\log(\varepsilon^{-1}))\sqrt{d}\sqrt{2^{m+1}}\\
  \leq & (M-L_m)(\lambda([0,y))+2^{-H})+(904.726+125.568\log(\varepsilon^{-1}))\sqrt{d}\sqrt{2^{m+1}}\\
  \leq & (M-L_m)\lambda([0,y))+(908.726+125.568\log(\varepsilon^{-1}))\sqrt{d}\sqrt{2^{m+1}}
 \end{aligned}
\end{multline}
on a set with probability at least $1-\varepsilon/6\cdot 2^{-2d}$.
Similarly by using (\ref{186}) instead of (\ref{187}) we obtain
\begin{multline}
 \label{212}
 \sum_{n=L_m+1}^M\mathbf{1}_{[0,y)}(x_n)\\
 \begin{aligned}
  \geq & \sum_{n=L_m+1}^M\mathbf{1}_{U_1}(q_n)\cdot \mathbf{1}_{V_1}(r_n)\\
  & +\sum_{h=1}^{H-1}\sum_{n=L_m+1}^M\left(\mathbf{1}_{U_{h+1}\backslash U_h}(q_n)\cdot \mathbf{1}_{V_{h+1}}(r_n)+\mathbf{1}_{U_h}(q_n)\cdot \mathbf{1}_{V_{h+1}\backslash V_h}(r_n)\right)\\
  \geq & (M-L_m)\lambda([0,y))-(908.726+125.568\log(\varepsilon^{-1}))\sqrt{d}\sqrt{2^{m+1}}
 \end{aligned}
\end{multline}
on the same set of probability bounded from below by $1-\varepsilon/6\cdot 2^{-2d}$.  Therefore we have proved (\ref{202}) which finally concludes the proof of the Theorem.





\vspace{3ex}

\small{DEPT. OF MATHEMATICS, BIELEFELD UNIV., P.O.Box 100131, 33501 Bielefeld, Germany\\
\textit{E-Mail address:} \url{tloebbe@math.uni-bielefeld.de}}


\begin{thebibliography}{99999}
\bibitem{AB10} Aistleitner, C.: Berkes, I.: On the central limit theorem for $f(n_kx)$, Probab. Theory Relat. Fields 146, 267-289 (2010)
\bibitem{A13} Aistleitner, C.: On the inverse of the discrepancy for infinite dimensional infinite sequences, J. Complexity 29, 182-194 (2013)
\bibitem{AFF13} Aistleitner, C., Fukuyama, K., Furuya, Y.: Optimal bound for the discrepancies of lacunary sequences, Acta Arith. 158, 229-243 (2013)
\bibitem{AW} Aistleitner, C., Weimar, M.: Probabilistic star discrepancy bounds for double infinite random matrices, 
Monte Carlo and Quasi-Monte Carlo Methods 2012, Springer (2013)
\bibitem{A04} Atanassov, E.I.: On the discrepancy of Halton sequences, Math. Balkanica, New Series 18, 15-32 (2004)
\bibitem{BS96} Bach, E., Shallit, J.: Algorithmic number theory. Vol. 1. Foundations of Computing Series. MIT Press, Cambridge, MA (1996)
\bibitem{CBR12} Conze, J.-P., Le Borgne, S., Roger, M.: Central limit theorem for stationary products of toral automorphisms, Discrete Contin. Dyn. Syst. 32, 1597-1626 (2012)
\bibitem{D07} Dick, J.: A note on the existence of sequences with small star discrepancy, J. Complexity 23, 649-652 (2007)
\bibitem{DP10} Dick, J.,Pillichshammer, F.: Digital nets and sequences, Discrepancy theory and quasi-Monte Carlo integration, Cambridge University Press, Cambridge, MA (2010)
\bibitem{DGKP08} Doerr, B., Gnewuch, M., Kritzer, P., Pillichshammer, F.: Component-by-component construction of low-discrepancy point sets of small size, Monte Carlo Methods Appl. 14, 129-149 (2008)
\bibitem{DT97} Drmota, M., Tichy, R.F.: Sequences, discrepancies and applications, vol. 1651 of Lecture Notes in Mathematics, Springer, Berlin, Heidelberg, New York (1997)
\bibitem{EG55} Erd\H{o}s, P., G\'al, I.S.:  On the law of iterated logarithm, Proc. Kon. Nederl. Akad. Wetensch. 58, 65-84 (1955)
\bibitem{EM96} Einmahl, U., Mason, D.M.: Some universal results on the behavior of increments of partial sums, Ann. Prob. 24, 1388-1407 (1996)
\bibitem{F08} Fukuyama, K.: The law of the iterated logarithm for the discrepancies of $\{\theta^nx\}$, Acta. Math. Hungar. 118, 155-170 (2008)
\bibitem{G66} Gaposhkin, V.F.: Lacunary series and independent functions, Russian Math. Surv. 21, 3-82 (1966)
\bibitem{G70} Gaposhkin, V.F.: The central limit theorem for some weakly dependent sequences, Theory Probab. Appl. 15, 649-666 (1970)
\bibitem{G08} Gnewuch, M.: Bracketing numbers for axis-parallel boxes and applications to geometric discrepancy, J. Complexity 24, 154-172 (2008)
\bibitem{H60} Halton, J.H.: On the efficiency of certain quasi-random sequences of points in evaluating multi-dimensional integrals, Numer. Math. 2, 84-90 (1960)
\bibitem{HNWW01} Heinrich, S., Novak E., Wasilkowski, G.W., Wo\'zniakowski, H.: The inverse of the star-discrepancy depends linearly on the dimension, Acta Arith. 96, 279-302 (2001)
\bibitem{H04} Hinrichs, A.: Covering numbers, Vapnik-\v Cervonenkis classes and bounds for the star-discrepancy, J. Complexity 20, 477-483 (2004)
\bibitem{K46} Kac, M.: On the distribution of values of sums of the type $\sum f(2^kt)$, Ann. Math. 47, 33-49 (1946)
\bibitem{K49} Kac, M.: Probability methods in some problems  of analysis and number theory, Bull. Am. Math. Soc. 55, 641-665 (1949)
\bibitem{L09} Lemieux, C.: Monte Carlo and quasi-Monte Carlo sampling, Springer (2009)
\bibitem{L13} Levin, M.: Central Limit theorem for $\mathbb{Z}^d_+$-actions by toral endomorphisms, Electron. J. Probab. 18, no. 35, 1-42 (2013)
\bibitem{N92} Niederreiter, H.: Random Number Generation and Quasi Monte-Carlo Methods, volume 63 of CBMS-NSF Regional Conference Series in Applied Mathematics, SIAM, Philadelphia, PA (1992)
\bibitem{P75} Philipp, W.: Limit theorems for lacunary series and uniform distribution mod 1, Acta Arith. 26, 241-251 (1975)
\bibitem{R80} Roth, K.F.: On irregularities of distribution I-IV, Mathematika 1, 73-79 (1954), Comm. Pure Appl. Math. 29, 739-744 (1976), Acta Arith. 35, 373-384 (1979) and Acta Arith. 37, 67-75 (1980)
\bibitem{SZ47} Salem, R., Zygmund, A.: On lacunary trigonometric series, Proc. Nat. Acad. Sci. USA 33, 333-338 (1947)
\bibitem{SZ50} Salem, R., Zygmund, A.: La loi du logarithme it\'er\'e pour les s\'eries trigo- nom\'etriques lacunaires, Bull. Sci. Math. 74, 209-224 (1950)
\bibitem{T61} Takahashi, S.: A gap sequence with gaps bigger than the Hadamards, Tohoku Math. J. 13, 105-111 (1961)
\bibitem{T62} Takahashi, S.: An asymptotic property of a gap sequence, Proc. Japan Acad. 38, 101-104 (1962)
\bibitem{WH00} Wang, X., Hickernell F.J.: Randomized Halton sequences, Math. Comput. Modelling 32, 887-899 (2000)
\bibitem{W59} Weiss, M.: The law of the iterated logarithm for lacunary trigonometric series, Trans. Amer. Math. Soc. 91, 444-469 (1959)
\bibitem{W16} Weyl, H.: \"Uber die Gleichverteilung von Zahlen mod. Eins, Math. Ann. 77, 313-352 (1916)
\bibitem{W55} Wigner, E.P.: Characteristic vectors of bordered matrices with infinite dimensions, Ann. of Math. 62, 548-564 (1955)
\end{thebibliography}
\end{document}